\documentclass[12pt]{amsart}
\usepackage{amsmath,amsthm,amsfonts,latexsym,amssymb,amscd,color,bm}
\usepackage{graphics,multimedia,tikz}
\usepackage{chemarr}
\newtheorem*{clm*}{Claim}
\theoremstyle{definition}

\numberwithin{equation}{section}
\newcommand{\cproof}{\noindent{\it Proof of Claim.}\ } 
\newcommand{\cqed}{\hfill\rule{1.3mm}{3mm}}

\newcommand{\wec}[1]{{\mathbf{#1}}}  

\newcommand{\al}[1]{{\mathbf{\uppercase{#1}}}}
\newcommand{\str}[1]{{\mathbb{\uppercase{#1}}}}

\DeclareMathOperator{\Pol}{Pol}
\DeclareMathOperator{\pPol}{pPol}

\DeclareMathOperator{\CSP}{CSP}

\newcommand{\var}[1]{\mathcal{#1}}
\DeclareMathOperator{\Clo}{Clo}

\newcommand{\pr}{\textrm{pr}}
\newcommand{\fin}{\textrm{fin}}

\newcommand{\comment}[1]{}

\newcommand{\OO}{\mathcal{O}}
\newcommand{\II}{\mathcal{I}}
\newcommand{\BB}{\mathcal{B}}
\newcommand{\CC}{\mathcal{C}}
\newcommand{\DD}{\mathcal{D}}
\newcommand{\MM}{\mathcal{M}}
\newcommand{\UU}{\mathcal{U}}

\newcommand{\ul}[1]{$\underline{{\text #1}}$}

\let\phi=\varphi
\let\epsilon=\varepsilon

\makeatletter
\@namedef{subjclassname@2020}{\textup{2020} Mathematics Subject Classification}
\makeatother

\begin{document}

\title[Maximal Clones, Minimal Clones]{Ivo G.~Rosenberg's Work on\\
  Maximal Clones and Minimal Clones}

\author{\'Agnes Szendrei}
\address[\'Agnes Szendrei]{Department of Mathematics\\
University of Colorado\\
Boulder, CO 80309-0395\\
USA}
\email{szendrei@colorado.edu}

\subjclass[2020]{Primary: 08A40; Secondary: 08B50, 08A05, 08A70, 68Q25}
\keywords{clone, completeness, maximal clone, Rosenberg's Completeness Theorem,
minimal clone, Rosenberg's Theorem on minimal clones}

\begin{abstract}
  The aim of this paper is to honor
  Ivo G.~Rosenberg by describing some of his most influential results
  and their impact in logic, discrete mathematics, algebra, and
  computer science.
\end{abstract}

\maketitle

\section{Introduction}
\label{sec-intro}

Ivo Rosenberg's published work touches on a wide range of topics in
mathematics, including logic and discrete mathematics, graph theory and
combinatorics, algebra, geometry, topology, and computer science.
A large portion of his results are concerned with \emph{clones}, which are
families $\CC$ of finitary operations (i.e., families of
functions in finitely many variables) on a fixed set $A$, such that 
$\CC$ is closed under composition and
contains all \emph{projection operations}
$\pr_i^{(n)}\colon A^n\to A$, $(a_1,\dots,a_n)\mapsto a_i$
($n\ge1$, $1\le i\le n$).
The largest clone on $A$ is the clone $\OO_A$ of all operations, and the
least clones on $A$ is the clone $\II_A$ of projections.

The study of clones entered into mathematics from
two independent sources: logic and algebra.
In the propositional calculus of $2$-valued logic, it has been known for
a long time that conjunction and negation
--- or similarly, implication and negation --- 
form a ``complete
set of logical connectives'' in the sense that the corresponding
Boolean functions (or truth functions),
together with the projections (which are the
`trivial' Boolean functions that record the truth values of the
variables),
are sufficient to generate all Boolean functions via composition.
Using the concept of a clone, we can restate the definition of completeness
as follows: a set $F$ of Boolean functions is \emph{complete} if $F$ generates
the clone of all Boolean functions, and 
a set of logical connectives is \emph{complete} if 
the set of the corresponding Boolean functions is complete.
In 1920, Post announced his classification of
all possible propositional calculi (or in modern terminology:
his classification of all clones of Boolean functions), which
implies the following
general completeness criterion
for $2$-valued logic (see \cite{Post-0, Post}):

\smallskip
\noindent
{\it
A set $F$ of Boolean functions is complete
if and only if there exist not necessarily distinct functions 
$f_1,\dots,f_5\in F$ such that $f_1$ is not monotone with respect to the order
$0\le1$, $f_2$ does not fix $0$, $f_3$ does not fix $1$, $f_4$ is not
invariant under switching the roles of $0$ and $1$, and $f_5$ is not linear
when expressed as a polynomial function over the $2$-element field
$\al{Z}_2$.}\footnote{Post used a slightly weaker notion than the concept of
a clone, which he called a \emph{closed set of functions}.
Every clone is a closed set of functions, but a closed set of functions may
not contain projections. However, if we add all projections to a closed set of
functions, we always get a clone.}

\smallskip

In algebra, the fundamental objects of study are
\emph{algebraic structures} (briefly: \emph{algebras}),
such as groups or rings.
Algebras are devices for computation, and can be formally described as
pairs $\al{A}=(A;F)$ where 
$A\not=\emptyset$ is the base set of $\al{A}$ and $F$ is the set of
basic operations of $\al{A}$.
Computations in $\al{A}$ that apply basic operations multiple times to
elements of $\al{A}$,
can be captured by \emph{derived operations} obtained
from the basic operations of $\al{A}$ and projection operations
by composition.
Hence the derived operations of $\al{A}=(A;F)$ form a 
clone, which is called \emph{the clone of $\al{A}$} and is
denoted by $\Clo(\al{A})$.
It is easy to see that $\Clo(\al{A})$ is the clone on
$A$ generated by $F$,
that is, $\Clo(\al{A})$ is the least clone on
$A$ containing $F$.
The concept of a clone was introduced by Philip Hall
(cf.~\cite[p.~127]{Cohn} and \cite[p.~126]{ALVIN3}),
who was apparently motivated by the difficulty of
doing calculations with derived operations of classical algebras,
especially groups and rings.

Rosenberg was a leader in the study of clones for over 50 years.
The aim of this paper is to honor
him by presenting two of his
most significant results --- one on maximal clones on finite sets
and the other on minimal clones ---, and discussing their profound impact
in logic, algebra, and computer science.

\section{Maximal Clones on Finite Sets and\\ Rosenberg's Completeness Theorem}
\label{sec-max-finite}

In this section we present
Ivo G.~Rosenberg's generalization of
Post's Completeness Theorem from $2$-valued
logic to $k$-valued logics for arbitrary finite $k\ge2$,
which is a landmark result in both logic and algebra.
He announced the theorem in 1965 in~\cite{Ros-max0}, and
the following year he submitted the result,
together with its 90-page proof, for publication, which appeared in
1970 in~\cite{Ros-max1}.

The concept of completeness carries over in a straightforward manner
from $2$-valued logic to $k$-valued logic for any $k\ge2$, as follows.
In $k$-valued logic ($k\ge2$) the \emph{truth functions} are the
finitary operations on
the $k$-element set $A$ of truth values, and a set $F$ of truth
functions --- that is, a set $F$ of operations on $A$ ---
is said to be \emph{complete}
if $F$ generates the clone $\OO_A$ of all operations on $A$.
A clone $\CC$ on a set $A$ is called \emph{maximal} if $\CC$ is a
proper subclone of $\OO_A$, and $\OO_A$ is the only clone on $A$ that
properly contains $\CC$.
If $k<2$, then every set of truth functions in $k$-valued logic is
complete, and if $|A|<2$, then maximal clones do not exist on $A$. 
Therefore, we will assume throughout that $k,|A|\ge2$, and in this section
we will also assume that $k$ and $A$ are finite.

\subsection{Background}\label{subsec-background}
To place Rosenberg's Completeness Theorem into context
let us review several important results on clones
that were discovered during the decades between
Post's work and the
announcement of Rosenberg's Theorem.
We start with two results that reveal essential differences
between $2$-valued logic and $k\,(>2)$-values logics.
In 1939, S{\l}upecki~\cite{Slup} proved the following
sufficient condition for completeness in $k\,(>2)$-valued logics, which
fails for $k=2$:

\smallskip
\noindent
{\it
If $A$ is a finite set of size $>2$ and $F\subseteq\OO_A$
contains all unary operations from $\OO_A$,
then $F$ is complete if and only if some operation in $F$ is both surjective
and not essentially unary (i.e., depends on at least two of its variables).}

\smallskip
\noindent
Another difference between $2$-valued and $k\,(>2)$-valued logics
that came to light in the late 1950's
is that the lattice of subclones of $\OO_A$ is 
much more complex for $|A|>2$ than for $|A|=2$.
Indeed, for $|A|=2$
Post proved in~\cite{Post} ---
by explicitly determining all subclones of $\OO_A$ ---
that $\OO_A$ has only countably many subclones.
However, for finite sets $A$ of size $>2$,
Yanov and Muchnik~\cite{YanMuc} proved that
$\OO_A$ has continuumly many subclones.

Nevertheless, other results that emerged in the 1950's 
gave hope that on the question of
completeness, the difference between
$2$-valued and $k\,(>2)$-valued logics may not be 
insurmountable. 
For example, Kuznetsov~\cite{Kuznetsov} proved that
for any finite set $A$ ($|A|\ge2$),
\begin{equation}\label{eq-kuzn}
  \begin{matrix}\text{\it every proper subclone of $\OO_A$ is contained
      in a maximal}\\
    \text{\it subclone of $\OO_A$, and the number of maximal subclones
      of $\OO_A$ is finite.}
  \end{matrix}  
\end{equation}  
Statement~\eqref{eq-kuzn}
implies that a set $F\subseteq\OO_A$ is complete if and only if $F$ is
not contained in any one of the finitely many maximal subclones of $\OO_A$.
Hence, finding an `efficient' completeness criterion for $k$-valued logics
($k\ge2$) is equivalent to finding
an `efficient' description for the maximal subclones of $\OO_A$ for all finite
sets $|A|\ge2$.
For the special case $|A|=2$, Post's Completeness Theorem is `efficient'
in this sense, and shows that $\OO_A$ has five maximal subclones.
For $|A|=3$, Yablonski\v{\i}~\cite{Yab}
presented a similar `efficient' completeness criterion by finding
all $18$ maximal subclones of $\OO_A$
(a result announced earlier without proof in~\cite{Yab-max}).

Parallel to these developments, but largely independently of them, 
in the early 1950's Foster initiated a study of those
algebras
$\al{A}=(A;F)$ where either the set $F$ of operations
itself, or $F$ together with all constant operations on $A$,
is a complete set of operations on $A$.
The nontrivial finite algebras\footnote{An algebra is \emph{nontrivial}
if it has more than one element.} $\al{A}$ of the first kind are called
\emph{primal}\footnote{%
In \cite{Foster-boolean1, Foster-boolean2} Foster called these algebras
\emph{functionally strictly complete}, and he
switched to the currently accepted terminology \emph{primal}
a few years later.}
algebras, the motivating examples being the $2$-element Boolean algebra and
finite fields of prime order.
The nontrivial finite algebras $\al{A}$
of the second kind
are called \emph{functionally complete} algebras, and they also include all 
finite fields.
Foster was primarily interested in the structure of all algebras
in the variety generated by a given primal or functionally complete
algebra $\al{A}$.\footnote{%
A \emph{variety} or \emph{equational class} is a class of algebras axiomatized
by a set of identities. The variety generated by $\al{A}$
is the class of all algebras that satisfy every identity true in $\al{A}$.}
For example, in \cite{Foster-boolean2} he proved that if $\al{A}$ is a primal
algebra, then every finite or infinite algebra in the variety
$\var{V}(\al{A})$ generated by
$\al{A}$ is a Boolean power\footnote{%
A Boolean power of $\al{A}$ is a special kind of subdirect product of copies
of $\al{A}$. 
In \cite{Foster-boolean1, Foster-boolean2} Foster called this construction
\emph{normal subdirect sum of $\al{A}$}.}
of $\al{A}$, and
hence $\var{V}(\al{A})$ is a \emph{minimal variety},
that is, it contains no other
proper subvariety than the trivial variety of one-element algebras.

The powerful, and nowadays ubiquitous, technique of using compatible relations
to describe clones on finite sets --- or, in a more algebraic
language, the technique of using the subalgebras of finite powers
of a finite algebra $\al{A}=(A;F)$ to determine the
clone of $\al{A}$ ---, which is based on the Galois connection between
operations and relations (see~\cite{BKKR, Gei}),
was not fully developed until after Rosenberg proved his completeness theorem. 
However, some related
terminology and notation will be useful throughout this article, therefore
we introduce them here.
If $A$ is a set, $f$ is an $n$-ary operation and $\rho$ is an $m$-ary relation
on $A$, we say that $f$ \emph{preserves} $\rho$
(or $\rho$ is \emph{invariant with respect to} $f$, or
$\rho$ is \emph{compatible with} $f$) if $\rho$ is closed under the operation
$f$ when $f$ is applied coordinatewise to elements of $\rho$.
Thus, saying that $f$ preserves $\rho$ is equivalent to saying that
$\rho$ is (the base set of) a subalgebra of the $m$-th power
algebra $(A;f)^m$, but 
it is also equivalent to saying that $f$ is an $n$-ary homomorphism
$(A;\rho)^n\to(A;\rho)$, or \emph{polymorphism}, of the relational
structure $(A;\rho)$.
It is easy to see that the set of all operations on $A$ that
preserve $\rho$ --- i.e., the set of all polymorphisms of $\rho$ ---
is a clone on $A$, which will be denoted by $\Pol(\rho)$.

Kuznetsov's proof of the fact \eqref{eq-kuzn} mentioned above yields
that for any finite set $A$ ($|A|\ge3$),
\begin{equation}\label{eq-kuzn2}
  \begin{matrix}\text{\it every maximal clone on $A$ is of the form
      $\Pol(\rho)$}\\
    \text{\it for some relation $\rho$ of arity $\le|A|$.}
  \end{matrix}  
\end{equation}

\subsection{Rosenberg's description of the maximal clones on finite sets}
The maximality
of $\Pol(\rho)$ was known before 1965  for the following familiar
relations on finite sets $A$ ($|A|\ge2$):
\begin{enumerate}
\item[(1)]
  $\rho$ is a partial order on $A$ with a least element and a greatest
  element~\cite{Mar};
\item[(2)]
  $\rho$ is (the graph of) a fixed point free permutation on $A$ such that
  all cycles of $\rho$ have the same prime length~\cite{Yab};
\item[(4)]
  $\rho$ is a nontrivial equivalence relation on $A$~\cite{Yab};\footnote{%
  The trivial equivalence
  relations on $A$ are the equality relation and $A\times A$;
  their polymorphism clones are $\OO_A$.}
\item[(5)$_1$]
  $\rho$ is a nontrivial unary relation on $A$~\cite{Yab};\footnote{%
  The trivial unary relations on $A$ are
  $\emptyset$ and $A$; their polymorphism clones are $\OO_A$.}
\item[(6)$_*$]
  $\rho$ is the relation
  \[
  \iota_A:=\{(a_1,\dots,a_{|A|})\in A^{|A|}: |\{a_1,\dots,a_{|A|}\}|<|A|\}
  \]
  where $|A|>2$~\cite{Yab}.
\end{enumerate}  

In particular, if $A$ is the set $\mathbf{2}:=\{0,1\}$
of truth values in $2$-valued logic,
then Items~(1), (2), and (5)$_1$ above yield
four maximal subclones in $\OO_A$, namely the clone
$\Pol(\le)$ of all Boolean functions that are monotone
with respect to the order $0\le 1$, the clone
$\Pol(\neg)$ of all Boolean functions that are invariant under
negation (i.e., under switching $0$ and $1$), and
for both truth values $v=0$ and $v=1$, the clone
$\Pol(\{v\})$ of all Boolean functions that fix $v$.
By Post's description of the maximal clones of Boolean functions,
there is a fifth maximal clone:
the clone of all Boolean functions that are linear, when expressed
as polynomial functions over the 2-element field $\al{Z}_2$.

Returning to the case when $A$ is any finite set ($|A|\ge2$), notice that
maximal clones of the types (4) and (6)$_*$ exist only if $|A|>2$.
Moreover, the maximality of the clone $\Pol(\iota_A)$
for the relation $\rho=\iota_A$ in (6)$_*$
is closely related to S{\l}upecki's Theorem.
Indeed,
it is easy to check that $\Pol(\iota_A)$ contains all operations on $A$
that are either essentially unary or non-surjective.
Therefore, by S{\l}upecki's Theorem, 
$\Pol(\iota_A)$ cannot contain any one of the remaining operations, and
hence $\Pol(\iota_A)$ is a maximal
subclone of $\OO_A$. 
For this reason, $\iota_A$ is often referred to as
\emph{S{\l}upecki's relation}, and $\Pol(\iota_A)$
as \emph{S{\l}upecki's clone} on $A$.

Rosenberg's celebrated completeness theorem expands the list of relations
(1)--(6)$_*$ above to a list of 
relations $\rho$ on finite sets $A$ such that the corresponding clones
$\Pol(\rho)$ yield the full list of
maximal clones on $A$. The relations come in six types.
Three of the of types, namely~(1), (2), and (4),
are exactly the relations described above.
The families~(5)$_1$ and (6)$_*$ above are enlarged considerably to get
Rosenberg's relations of types~(5) and (6), and relations of type~(3)
are added to the list. This last type is not entirely new either, because
the maximal clone of linear Boolean functions
in Post's Theorem turns out to have
the form $\Pol(\rho)$ for the unique relation $\rho$ of type~(3)
on $A=\mathbf{2}$.

\medskip
\noindent
{\bf Rosenberg's Completeness Theorem {\rm\cite{Ros-max0, Ros-max1}}.}
{\it
  Let $A$ be a finite set $(|A|\ge3)$.
  The maximal clones on $A$ are exactly
  the clones $\Pol(\rho)$ where $\rho$ is one of the relations described in
  $(1),(2),(4)$ above or $(3),(5),(6)$ below:
  \begin{enumerate}
  \item[(3)]
  for some elementary abelian $p$-group $(A;+)$,
  $\rho$ is the $4$-ary relation
  \begin{equation}\label{eq-affine-rel}
    \alpha_+:=\{(a_1,a_2,a_3,a_4)\in A^4:a_1-a_2+a_3=a_4\},
  \end{equation}  
  which is the graph of the `ternary addition'
  $x_1-x_2+x_3$; 
  \item[(5)]
  for some integer $m$ $(1\le m<|A|)$, $\rho$ is an $m$-ary
  totally reflexive, totally symmetric\footnote{%
  An $m$-ary relation $\rho$ on $A$ is \emph{totally reflexive}
  if it contains all tuples $(a_1,\dots,a_m)\in A^m$ where $a_1,\dots,a_m$
  are not pairwise distinct, and \emph{totally symmetric} if $\rho$ is
  invariant under all permutations of its coordinates.}
  relation on $A$ such that
  $\rho\not=A^n$ and for at least one element $c\in A$ we have that
  $\{c\}\times A^{m-1}\subseteq\rho$;
  \item[(6)]
    for some integers $m$ and $r$ $(3\le m\le|A|,\ r\ge1)$,
    and some $r$-element set
    $T=\{\theta_1,\dots,\theta_r\}$ of
    equivalence relations on $A$ such that
  \begin{itemize}
  \item
    each $\theta_i$ $(1\le i\le r)$ has exactly $m$ blocks, and
  \item
    $\bigcap_{i=1}^r B_i$ is nonempty whenever $B_i$ is a block of
    $\theta_i$ for every $i$,
  \end{itemize}
  $\rho$ is the relation $\lambda_T$ defined as follows:
  \begin{multline}\label{eq-mregular-rel}
    \lambda_T:=\{(a_1,\dots,a_m)\in A^m: \text{for every $i$ $(1\le i\le r)$,
      there exist}\\
    \text{two elements among $a_1,\dots,a_m$ that are in
      the same block of $\theta_i$}\}.
  \end{multline}
  \end{enumerate}
  Consequently, a set $F$ of functions on $A$ is complete --- or equivalently,
  an algebra $(A;F)$ is primal --- if and only if $F\not\subseteq\Pol(\rho)$
  holds for every relation $\rho$ in $(1)$--$(6)$.
}    

\medskip

The relations $\rho$ in (3), (5), and (6) are called
\emph{affine relations}, \emph{central relations}, and
\emph{$m$-regular relations}\footnote{%
An alternative name in the literature for $m$-regular relations is
\emph{$m$-universal relations}.},
respectively. It is easy to see that affine relations exist on $A$ if and
only if $|A|$ is a prime power.
Observe also that the families of relations in (5) and (6) do indeed extend
the families in (5)$_1$ and (6)$_{*}$, as claimed earlier. Indeed, 
the relations in (5)$_1$ are exactly the unary central relations
on $A$; this is what the subscript~${}_1$ appended to (5)
is intended to indicate.
Furthermore, S{\l}upecki's relation $\iota_{|A|}$
in (6)$_*$ is the unique $|A|$-regular relation on $A$, namely
$\lambda_{\{=\}}$;
here the subscript~${}_*$ appended to (6) is intended to indicate that
S{\l}upecki's relation is a distinguished member of the set of relations
in (6).
We will say more about the affine, central, and
$m$-regular relations and the maximal clones they determine
in Subsection~\ref{subsec-356}.

To accompany his completeness theorem above, 
Rosenberg proved in \cite{Ros-max2}
that almost all maximal clones $\Pol(\rho)$ where $\rho$ is one of
the relations in (1)--(6) are distinct.
More precisely, he proved that
if $\rho$ and $\sigma$ are distinct relations in (1)--(6),
then $\Pol(\rho)=\Pol(\sigma)$ holds if and only if either
(i)~$\rho$ and $\sigma$ are both of type (1) and they are inverses
of each other, or else
(ii)~$\rho$ and $\sigma$ are both of type (2) and they are powers of each other.
In \cite{Ros-max3},
Rosenberg also found a
formula for computing the number of maximal clones on an $n$-element set.

\subsection{Idea of the proof of Rosenberg's Completeness
  Theorem}\label{subsec-ideaofproof}
No easy proof is known for Rosenberg's description of all maximal
clones on finite sets.
Here, we will only discuss the idea and the main
difficulties of the original proof in \cite{Ros-max1}.\footnote{Detailed
proofs for Rosenberg's Completeness Theorem are also presented in
\cite{Quackenbush}, \cite{Lau-book}, and \cite{ALVIN3}.}
The starting point of the proof is the fact in \eqref{eq-kuzn2} that  
every maximal clone on a finite
set $A$ ($|A|\ge3$) is a member of the finite family
$\{\Pol(\rho):\rho\in\mathcal{R}\}$
of clones where $\mathcal{R}$
is the set of relations $\rho$ of arity $\le|A|$ on $A$ such that
$\Pol(\rho)\not=\OO_A$ (i.e., $\rho\not=\emptyset$ and
$\rho$ is not a diagonal relation\footnote{%
An $m$-ary relation $\rho$ on $A$
is \emph{diagonal} if there exists an equivalence
relation $\epsilon$ on the set $\{1,\dots,m\}$
such that for all $(a_1,\dots,a_m)\in A^m$,
we have $(a_1,\dots,a_m)\in\rho$
if and only if $a_i=a_j$ holds whenever $i\,\epsilon\,j$.
In particular, $A^m$ is an $m$-ary diagonal relation on $A$ that we get when
$\epsilon$ is the equality relation.
}).
The proof is an intricate elimination process, removing relations from
$\mathcal{R}$ step-by-step in such a way that
\begin{enumerate}
\item[--]
a relation $\sigma\in\mathcal{R}$ is eliminated only if there remains
a relation $\rho\in\mathcal{R}$ with $\Pol(\rho)\supseteq\Pol(\sigma)$  
still to be processed --- hence, no maximal clones are lost in the process ---,
and
\item[--]
at the end, for all relations $\rho$ that have not been eliminated,
$\Pol(\rho)$ is a maximal clone on $A$.
\end{enumerate}
Throughout the proof,
unless stated otherwise, when a clone is written in the form $\Pol(\rho)$,
$\rho$ is chosen to have minimum arity among the relations determining
$\Pol(\rho)$.

The proof has two main phases.

In Phase~1, Rosenberg isolates the maximal clones of
types (5)$_1$, (1), (2), and (3), and proves that
every maximal clone on $A$ other than those of types (5)$_1$, (1), (2), and (3) 
are of the form $\Pol(\tau)$ for a non-diagonal
totally reflexive, totally symmetric
relation $\tau$ of arity $\ge2$.
The difficulty in this part of the elimination process comes when Rosenberg
narrows down the possibilities for $\Pol(\rho)$ to the case where
$\rho$ has arity $\ge3$, 
$\Pol(\rho)$ is not contained in any of the maximal clones of types
(5)$_1$, (1), or (2), and $\Pol(\rho)$ is not contained in $\Pol(\tau)$ for any
non-diagonal totally reflexive, totally symmetric relation $\tau$ of
arity $\ge2$.
He proves that such a clone $\Pol(\rho)$
has to be maximal of type (3), although for
$p>2$ some ternary relations --- such as the graph of the 
binary operation $2x_1-x_2=x_1-x_2+x_1$ ---
also determine the corresponding affine maximal clone $\Pol(\alpha_+)$
in (3).

In Phase 2, it remains to consider the clones $\Pol(\tau)$
($\tau\in\mathcal{R}$)
where $\tau$ is a non-diagonal totally reflexive, totally symmetric
relation of arity $\ge2$ on $A$. The maximal clones of the types (4) and (5)
are isolated early on in this phase, and it is also established that
for all remaining maximal clones $\Pol(\tau)$, where $\tau$ is a
non-diagonal totally reflexive, totally symmetric relation, one can choose
$\tau$ to be a homogeneous\footnote{%
An $m$-ary relation $\tau$ on $A$ is called \emph{homogeneous} if
$(a_1,\dots,a_m)\in\tau$ holds whenever $a_1,\dots,a_m\in A$ and
there exists $u\in A$ such that
$(a_1,\dots,a_{i-1},u,a_{i+1},\dots,a_m)\in\tau$ for all $1\le i\le m$.}
relation.
The rest of the elimination process, which removes all homogeneous relations
other than those of type (6) from consideration,
takes about 25 pages, and is a real tour de force.

\subsection{More about the maximal clones of types (3), (5), and
  (6)}\label{subsec-356}
Functions and operations that are monotone with respect to a
partial order, commute with a permutation, preserve an equivalence
relation or preserve a subset occur in almost every area of mathematics,
therefore the maximal clones of types (1), (2), (4), and (5)$_1$ feel natural
and familiar. So are the maximal clones of type~(3) if one uses
the following algebraic description of them from \cite{Ros-max1}:
If $\al{V}:=(\al{Z}_p)^d$ is the $d$-dimensional vector space over the
$p$-element field $\al{Z}_p$ ($p$ prime),
$(V;+)$ is the additive group of $\al{V}$,
and $\alpha_+$ is the corresponding affine relation, then 
the members of the maximal clone $\Pol(\alpha_+)$
are exactly the operations of the form
\begin{equation}\label{eq-affine-max-clones}
\sum_{i=1}^n M_ix_i+v\quad (n\ge1),
\end{equation}  
where $v\in V$ and 
the coefficients $M_1,\dots,M_n$ are $d\times d$ matrices over $\al{Z}_p$.
Any other maximal clone of type~(3) associated to an elementary abelian
$p$-group $(A;+)$ of order $p^d$ is obtained from this one by renaming the
elements of the base set via an isomorphism $(V;+)\to(A;+)$.

The remaining relations in Rosenberg's list, i.e., 
the non-unary central relations and the
$m$-regular relations, are  certainly less familiar, and when encountered
for the first time, might look
mysterious.
Our goal here is to build some intuition for them and
the maximal clones they determine.
A common property of all maximal clones 
determined by a central relation $\rho$ of
arity $m\ge2$ or by an $m$-regular relation $\lambda_T$ ($m\ge3$)
on a finite set $A$ is that they contain
all operations on $A$ whose range has size $<m$, because these relations are
totally reflexive. In addition, every maximal clone determined by a central
relation $\rho$ of arity $m$ also contains all operations on $A$ whose
range has size $m$ and has an element from the \emph{center}
\begin{equation*}
  C_\rho:=\{c\in A:\{c\}\times A^{m-1}\subseteq\rho\}
\end{equation*}  
of $\rho$.

On the other hand, if
one focuses on surjective operations, the maximal clones
determined by $m$-ary central relation ($m\ge2$)
and the maximal clones determined by $m$-regular relations ($m\ge3$)
behave quite differently.
For central relations $\rho$, every surjective operation in $\Pol(\rho)$
also preserves $C_\rho$, which is a unary relation in (5)$_1$.
Therefore, in this context, central relations
may be thought of as higher arity variants
of the unary relations in (5)$_1$.
Note, however, that for surjective operations, 
preserving $\rho$ is a stronger requirement than preserving $C_\rho$.

Before discussing surjective operations in the maximal clones
$\Pol(\lambda_T)$ determined by $m$-regular relations, we will present 
an alternative description for these relations\footnote{%
This is the description that was
used by Rosenberg in \cite{Ros-max1}.},
which shows more explicitly
how the $m$-regular relations $\lambda_T$ are related to
S{\l}upecki's relation $\iota_B$ on an $m$-element set $B$.
First, notice that the two requirements on $T$ itemized in (6) are
equivalent to requiring that there is a surjective function
$\phi\colon A\to B^r$ mapping $A$ onto the product of $r$ copies of $B$
in such a way that $\theta_1,\dots,\theta_r$
are the kernels of the $r$ functions
$A\stackrel{\phi}{\to}B^r\stackrel{\pr_i^{(r)}}{\to}B$ ($1\le i\le r$)
where $\pr_i^{(r)}$
is the $i$-th $r$-ary projection operation on $B$.
Now, it is not hard to check that if $(\iota_B)^r$ denotes
the relation on $B^r$ which is obtained by
applying $\iota_B$ coordinatewise on $B^r$, then
$\lambda_T$ is the inverse image of $(\iota_B)^r$ under $\phi$.

For simplicity, we will discuss the surjective operations in $\Pol(\lambda_T)$
only in the case where $\phi$ is a bijection.
By renaming the elements of $A$
via the bijection $\phi$, we may assume without loss of generality
that $A=B^r$ and $\rho=(\iota_B)^r$. 
Now, using the fact that $\Pol(\iota_B)$ is S{\l}upecki's clone, and hence
all surjective operations in it are essentially unary, one can derive that
every surjective operation\footnote{%
In fact, the requirement for the operation to be \emph{surjective}
may be weakened to \emph{surjective in every coordinate of $B^r$}.}
in $\Pol\bigl((\iota_B)^r\bigr)$
is a `selector function' of the form
\begin{equation}\label{eq-reg-max-clones}
  \left(
  \begin{bmatrix} b_{11}\\ \vdots \\ b_{r1}\end{bmatrix},\ \dots,\
  \begin{bmatrix}b_{1n}\\ \vdots \\ b_{rn}\end{bmatrix}
  \right)\quad\mapsto\quad
  \begin{bmatrix}g_1(b_{i_1 j_1})\\ \vdots \\ g_r(b_{i_r j_r})\end{bmatrix}
  \qquad (n\ge1),
\end{equation}  
where the elements of $B^r$ are written as column vectors, $g_1,\dots,g_r$
are unary operations on $B$,
and the subscript $(i_\ell,j_\ell)$ describes for each
$1\le \ell\le r$ which coordinate of which of the $n$ arguments
of the operation is used to
compute the $\ell$-th argument of the result.

\section{The Impact of Rosenberg's Completeness Theorem}\label{sec-impact}

Over the past $50^+$ years, Rosenberg's theorem has had a
vast impact on research within the broader topic of `completeness'
and also far beyond.
Our goal here is to highlight the influence of the theorem and its proof
in logic, discrete mathematics, algebra, and computer science.

\subsection{Completeness criteria for special sets of functions}
One of the first applications of Rosenberg's Completeness Theorem was
Rousseau's characterization of Sheffer functions on finite sets. An operation
$f$ on a set $A$ ($|A|\ge2$)
is called a \emph{Sheffer function} if the singleton set
$\{f\}$ is complete, that is, the algebra $(A;f)$ whose only operation is $f$
is primal. In the special case
$A=\mathbf{2}$, the Boolean functions NAND and NOR are well-known examples of
Sheffer functions. Rousseau observed that Rosenberg's theorem
simplifies considerably when applied to one-element
sets of operations, and proved the following theorem in \cite{Rous}:

\smallskip
\noindent
{\it A function $f$ on a finite set
$A$ $(|A|\ge2)$ is a Sheffer function if and only if 
$f\notin\Pol(\rho)$ holds for the Rosenberg relations $\rho$ in
(2), (4), and (5)$_1$.}

\smallskip
\noindent
Equivalently:
A finite algebra $\al{A}=(A;f)$ ($|A|\ge2$) with a single
operation $f$ is primal if and only if $\al{A}$ has no nontrivial
automorphisms, no proper subalgebras, and no nontrivial congruences.

It is a straightforward consequence of Rosenberg's theorem that
a finite algebra $(A;F)$ ($|A|\ge2$)
is functionally complete if and only if $F\not\subseteq\Pol(\rho)$
holds for all Rosenberg relations $\rho$ in (1), (3), (4), (6), and all
non-unary relations $\rho$ in (5). In ~\cite{Ros-func-compl-surj-algs},
Rosenberg gave a necessary and sufficient condition --- along the lines
of Rousseau's result --- for the functional completeness of finite algebras
$(A;F)$ where $F$ contains only one operation or every operation
in $F$ is surjective.

\subsection{Completeness in algebra 1: Quasi-primal algebras and
  the Baker--Pixley Theorem}\label{subsec-quasiprimal}
Until the late 1970's, there seems to have been little interaction between
the research community studying completeness from the viewpoint of
multiple-valued logic and discrete mathematics,
and the researchers in general algebra who followed in
Foster's footsteps to investigate completeness. In general algebra the
emphasis was on the varieties generated by finite algebras that satisfy
one of several weaker conditions of completeness than primality, and
the approach was greatly influenced by Mal'tsev conditions, which
describe structural properties of varieties by the
existence of certain operations in the clone of a generating algebra.
The most successful generalization of primality turned out to be the
notion of quasi-primality, introduced by Pixley~\cite{Pixley-quasiprimal}.
An algebra $\al{A}=(A;F)$ is called \emph{quasi-primal} if
$\al{A}$ is finite, nontrivial (i.e., $|A|\ge2$), and
for every operation $f$ on $A$, $f$ belongs to
the clone of $\al{A}$ if and only if $f$ preserves all isomorphisms
between (not necessarily distinct) subalgebras of $\al{A}$.\footnote{%
In \cite{Pixley-quasiprimal}, Pixley called these algebras
\emph{simple algebraic algebras}, and introduced 
the current terminology \emph{quasi-primal algebras} for them
in \cite{Pixley-discriminator}.
Note also that the definition of quasi-primal algebras in these papers is
slightly different from the one given here, but equivalent.}
In the papers \cite{Pixley-quasiprimal} and \cite{Pixley-discriminator},
Pixley presented the following characterizations of quasi-primal algebras:

\smallskip
\noindent
{\it
A nontrivial finite algebra $\al{A}$ is quasi-primal
\begin{enumerate}
\item[$\Leftrightarrow$]
  every nontrivial subalgebra of $\al{A}$ is simple\footnote{%
  An algebra is said to be \emph{simple} if it 
  is nontrivial and has no nontrivial congruences.},
  and the clone of $\al{A}$
  contains ternary operations $p$ and $m$ such that the identities
  \begin{equation}\label{eq-maltsev}
    p(x,x,y)=p(y,x,x)=y,
  \end{equation}
  \begin{equation}\label{eq-majority}
    m(x,x,y)=m(x,y,x)=m(y,x,x)=x
  \end{equation}  
  hold in $\al{A}$;
\item[$\Leftrightarrow$]
  every nontrivial subalgebra of $\al{A}$ is simple, and the clone of $\al{A}$
  contains a ternary operation $q$ such that the identities
  \begin{equation}\label{eq-arithm}
    q(x,x,y)=q(y,x,y)=q(y,x,x)=y
  \end{equation}  
  hold in $\al{A}$;
\item[$\Leftrightarrow$]
  the clone of $\al{A}$ contains the ternary discriminator function $t$ on $A$
  defined by $t(a,b,c)=c$ if $a=b$ and $t(a,b,c)=a$ otherwise $(a,b,c\in
  A)$.
\end{enumerate} } 

\smallskip
\noindent
Clearly, the discriminator function $t$ satisfies the identities for $q$ in
\eqref{eq-arithm}.

An operation $p$ satisfying the identities in
\eqref{eq-maltsev} is called a \emph{Mal'tsev operation}, and by
\cite{Maltsev},
the existence of a Mal'tsev operation in the clone of an algebra $\al{A}$
is equivalent to the condition that congruences of algebras in the variety
$\var{V}(\al{A})$ generated by $\al{A}$ permute.
An operation $m$ satisfying the identities in
\eqref{eq-majority} is called a \emph{majority operation}, and by
results in 
\cite{Pixley-arithmetical} and \cite{Pixley-quasiprimal},
the existence of 
Mal'tsev as well as majority operations in the clone of an algebra $\al{A}$,
and the existence of an operation $q$ as in \eqref{eq-arithm}
in the clone of an algebra $\al{A}$, are both
equivalent to the condition that congruences of algebras in the variety
$\var{V}(\al{A})$ generated by $\al{A}$ permute and obey the distributive law.

Higher arity versions of majority operations are operations $u$ of any arity
$n\ge3$ which satisfy the identities
\begin{equation}\label{eq-nu}
  u(y,x,\dots,x)=\dots=u(x,\dots,x,y,x,\dots,x)=\dots=
  u(x,\dots,x,y)=x
\end{equation}
for every position of the `lone dissenter' $y$.
They are called \emph{near unanimity operations}, and they gained prominence
due to the 1975 paper \cite{Baker-Pixley} of Baker and Pixley, from which
we cite two results. The first one reveals the structural property of
varieties encoded by the existence of a near unanimity operation. 
Namely, the existence of a $(d+1)$-ary near unanimity operation in the
clone of an algebra $\al{A}$ is equivalent to the condition that
the following version of the Chinese Remainder Theorem
holds in the variety $\var{V}(\al{A})$ generated by $\al{A}$:
in any algebra in $\var{V}(\al{A})$, if a finite  system
$x_i\equiv a_i\ \text{\rm mod}\,\Theta_i$ ($1\le i\le r$, $r\ge d$)
of congruences is solvable $d$ at a time, then the whole system
is solvable. (This result is attributed to Huhn, cf.\
\cite[p.~169]{Baker-Pixley} and \cite[p.~89]{Huhn-weakdistr2}.)
The second result is about individual algebras, and is
usually referred to as the Baker--Pixley Theorem:

\smallskip
\noindent
{\it
  If $\al{A}=(A;F)$ is a finite algebra
  and its clone contains a $(d+1)$-ary near unanimity
  operation, then an operation $f$ on $A$ is in the clone of $\al{A}$ if and
  only if every subalgebra of $\al{A}^d$ is closed under $f$ (or equivalently,
  $f$ preserves all $d$-ary compatible relations of $\al{A}$).}

\smallskip
\noindent  
Two notable immediate consequences of the Baker--Pixley Theorem
are that (a)~for every fixed near unanimity
operation $u$ on a finite set there are only finitely many clones on
$A$ that contain $u$, and hence (b)~if a clone 
on a finite set contains a near unanimity operation, then it is
finitely generated.\footnote{%
A clone $\CC$ on $A$ is \emph{finitely generated} if for some finite set
$F$ of operations on $A$, $\CC$ is the smallest clone containing $F$.}

\subsection{Completeness in algebra 2: Abelianness and para-primal
  algebras}\label{subsec-abelian}
Another significant development in general algebra in the mid to late 1970's,
which is relevant to completeness, was the emergence of commutator theory,
and in particular, the isolation of the notion of abelianness
for arbitrary algebras and congruences. 
The theory generalizes classical commutator theories (e.g., the group
theoretical commutator and ideal multiplication in ring theory), and was 
first developed by Smith~\cite{Smith} for varieties which have Mal'tsev
term operations. Within a few years, the theory was
extended by Hagemann--Herrmann~\cite{Hagemann-Herrmann},
Gumm~\cite{Gumm-easyway}
and Freese--McKenzie~\cite{Freese-McKenzie-rs} to varieties
where every algebra has a modular congruence lattice, and later far beyond.
Typical examples of abelian algebras are (i)~modules over rings
(abelian groups in group theory and zero rings in ring theory),
and more generally,
(ii)~algebras whose clone is a subclone of the clone of the constant expansion
of a module, and also (iii)~all subalgebras of the algebras in (ii).
If the clone of such an algebra also contains a Mal'tsev
operation (the only choice being the ternary addition $x_1-x_2+x_3$ of
the abelian group of the module), then the algebra is called an
\emph{affine algebra}.
Every maximal clone of type~(3) in Rosenberg's Theorem is
the clone of an affine algebra. Indeed, with the notation used in
Subsection~\ref{subsec-356}, the module corresponding to this algebra
is obtained by considering
the vector space $\al{V}=(\al{Z}_p)^d$ over $\al{Z}_p$ as a module over
the full $d\times d$ matrix ring with entries in $\al{Z}_p$.
The clone of this module consists of all operations in
\eqref{eq-affine-max-clones} with $v=0$, and
it is not hard to see that by adding
all constant operations to this module we get an algebra
whose clone consists exactly of
the operations in \eqref{eq-affine-max-clones}.
Clearly, both clones include $x_1-x_2+x_3$.

Early on in the development of commutator theory,
Clark and Krauss~\cite{Clark-Krauss1} introduced a new
completeness notion for algebras, which is a common generalization of
quasi-primal algebras and abelian groups of prime order. They defined
an algebra $\al{A}=(A;F)$ to be \emph{para-primal} if $\al{A}$ is
a nontrivial finite algebra such that for every subalgebra $\al{B}$
of a finite power $\al{A}^I$ of $\al{A}$ and for every minimal subset
$J\subseteq I$ of coordinates of $\al{B}$ with the property that
the map $\pr_J$ projecting $\al{B}$ onto its coordinates in $J$
is one-to-one, we have that the image $\pr_J(\al{B})$ is the full
direct product $\prod_{j\in J}\pr_j(\al{B})$.
One of the main results of \cite{Clark-Krauss1}
extends Pixley's first two characterizations of quasi-primal algebras
mentioned above as follows:

\smallskip
\noindent
{\it
A nontrivial finite algebra $\al{A}$ is para-primal
\begin{enumerate}
\item[$\Leftrightarrow$]
  every nontrivial subalgebra of $\al{A}$ is simple, and the clone of $\al{A}$
  contains a Mal'tsev operation.
\end{enumerate}}

\smallskip
\noindent
Simultaneously, while investigating minimal varieties
with a Mal'tsev term operation,
McKenzie proved in \cite{McKenzie-ms1976} the following generalization
and strengthening of the theorem by Maurer and Rhodes~\cite{Maurer-Rhodes}
saying that every finite simple group is either functionally complete or
abelian:

\smallskip
\noindent
{\it
If $\al{A}$ is a nontrivial finite algebra such that $\al{A}$ is
\emph{strictly simple} (that is, $\al{A}$ is simple and has no nontrivial
proper subalgebras) and its clone contains a Mal'tsev operation, then
$\al{A}$ is either quasi-primal or affine.}

\smallskip
\noindent
A combination of these results, with further work, has revealed that,
except for having some strictly simple affine subalgebras,
para-primal algebras are very much like quasi-primal algebras
(see \cite{McKenzie-paraprimal} and \cite[Ch.~4]{szendreiBOOK}).
For example, the following is proved in \cite[p.~103]{szendreiBOOK}:

\smallskip
\noindent
{\it
If $\al{A}$ is a para-primal algebra and its affine subalgebras
  are $\al{B}_1,\dots,\al{B}_r$, then the clone of $\al{A}$
  contains a ternary operation $f$ such that
  the restriction of $f$ to each affine subalgebra $\al{B}_i$
  is the unique Mal'tsev operation in the clone of $\al{B}_i$,
  and for all other triples $(a,b,c)\in A^3$,
  the value of $f$ is computed by the ternary discriminator
  function $t$ defined earlier.}

\smallskip
\noindent
In \cite{Quackenbush}, Quackenbush used these ideas to give a new proof
for Rosenberg's Completeness Theorem. This approach simplified
many of the original arguments of Rosenberg, especially those in Phase~1
(see Subsection~\ref{subsec-ideaofproof}).

\subsection{Properties of maximal clones}
A lot of questions about maximal clones have been studied
since Rosenberg completed his description of all maximal clones on
finite sets.
Here we will only mention a few of them. 
Following up on the result of Yanov and Muchnik on the number of all clones on
finite sets $A$ with $|A|\ge3$, it was natural to wonder how this result
can be refined by counting the clones below each maximal clone one-by-one.
It is not hard to show that maximal clones of type~(3) have only
countably many subclones, 
finitely many if $|A|$ is prime (Salomaa~\cite{Salomaa}) and
infinitely many if $|A|$ is a prime power that is not prime.
On the other hand, all maximal clones of the remaining five types have
continuumly many subclones.
For the maximal clones of type~(2) the proof is much more difficult
than for the other types, and was found by
Demetrovics and Hann\'ak~\cite{Demetrovics-Hannak} and
by Marchenkov~\cite{Marchenkov}.

The proof of \eqref{eq-kuzn}  (see \cite{Kuznetsov} and
\cite{Ros-max1}) is based on the observation
that for finite sets
$A$ with $|A|\ge3$, the maximal clones $\MM$ on $A$
are determined by their unary parts $\MM^{(1)}$
consisting of all unary operations in $\MM$. Clearly, the unary part
of S{\l}upecki's clone, which is the set of all unary operations on $A$,
properly contains the unary parts of all the other maximal clones, but are
there any other proper inclusions among the unary parts of maximal clones?
This question
was investigated by Ma\v{s}ulovi\'c and Pech (n\'ee Ponjavi\'c) in
\cite{Ponj-Masul}, \cite{Ponj},
\cite{Masul-Pech}, \cite{Pech-Masul}, and \cite{Pech}.
Although some open problems still remain, especially when at least one of
the relations $\rho,\sigma$ is non-unary of type~(5) or (6),
the results they obtained show that proper inclusions
$\Pol^{(1)}(\rho)\subsetneq\Pol^{(1)}(\sigma)$ among the unary parts of
maximal clones $\Pol(\rho)$ and $\Pol(\sigma)$ are rare. For example,
$\sigma$ cannot be of type~(1) or (2) in such an inclusion, and it can be
of type~(3) only if $\rho$ is of type~(2) and either $|A|$ is prime or
$|A|=4$.
Surprisingly, they also found that there exist long chains of proper
inclusions among the unary parts of
maximal clones where the relations are
either all of type~(5) or all of type~(6). The lengths of the chains are
$|A|-1$ and $O(\sqrt{|A|})$, respectively.%

Which maximal clones $\MM$ on a finite set $A$ ($|A|\ge3$)
are finitely generated?
This question is significant, because it determines whether or not 
it is feasible to expect a Rosenberg-type description for the maximal
subclones of $\MM$. Indeed, if $\CC$ is a finitely generated clone on $A$,
then 
the idea of proof for Statements \eqref{eq-kuzn}--\eqref{eq-kuzn2}
mentioned in the preceding paragraph extends\footnote{%
The key is to replace the role of the unary part $\CC^{(1)}$ of the clone
by the $m$-ary part $\CC^{(m)}$, where $m$ is the maximum arity of the
operations in a generating set of $\CC$.} 
from $\OO_A$ to $\CC$ to prove that
\begin{itemize}
\item
  every proper subclone of $\CC$ is contained in a maximal subclone of
  $\CC$,
\item
  the number of maximal subclones of $\CC$ is finite, and
\item
  every maximal subclone of $\CC$ is of the form $\CC\cap\Pol(\rho)$
  for some relation $\rho$ on $A$.
\end{itemize}
Conversely, it is easy to see that
the first two conditions here imply that
$\CC$ is finitely generated.

It turns out that 
most maximal clones on $A$ are finitely generated, but not all, unless the set
$A$ is small. In more detail, all maximal clones
$\MM=\Pol(\rho)$ on $A$ ($|A|\ge3$)
where $\rho$ is $\underline{\text{not}}$ of type~(1)
are finitely generated. For $\rho$ of type~(3)
this follows from Rosenberg's description
of these clones in \eqref{eq-affine-max-clones}, and for the remaining types
this follows
from the results of Schofield~\cite{Schofield} and Lau~\cite{Lau-fingen}.
The maximal clones $\MM=\Pol(\rho)$ where $\rho$ is of type (2), (4), or (5)
contain near unanimity operations, therefore the Baker--Pixley Theorem
also implies that these maximal clones are finitely generated.
For the remaining case of type~(1), that is,
for the maximal clones of the form
$\MM=\Pol(\le)$ where $\le$ is a bounded partial order on $A$,
the problem asking
which of them are finitely generated is still open, in general.
A large class of finite bounded posets $(A;\le)$ for which
$\MM=\Pol(\le)$ is known to be finitely generated are 
the posets that can be obtained from finite lattices
by removing a (possibly empty) convex subset. For these posets
Demetrovics, Hann\'ak, and R\'onyai~\cite{Demetrovics-Hannak-Ronyai}
proved that
$\Pol(\le)$ contains a near unanimity operation, so the Baker--Pixley
Theorem applies. This result implies, in particular, that $\Pol(\le)$ is
finitely generated for every bounded partial order on $A$,
provided $|A|\le7$.
Tardos~\cite{Tardos} found the first example 
of a bounded partial order $\le$ on a finite set for which $\Pol(\le)$
is not finitely generated; it is the order depicted in Figure~1
on an $8$-element set.

\setlength{\unitlength}{1truemm}
  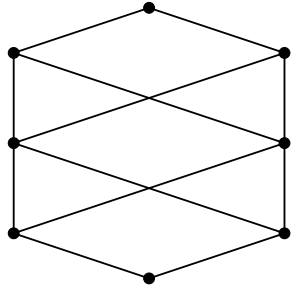
\begin{figure}
    \begin{picture}(40,40)
      \put(0,0){
\begin{tikzpicture}[scale=.12]
  \draw[fill] (0,5) circle (0.6);
  \draw[fill] (-15,10) circle (0.6);
  \draw[fill] (15,10) circle (0.6);
  \draw[fill] (-15,20) circle (0.6);
  \draw[fill] (15,20) circle (0.6);
  \draw[fill] (-15,30) circle (0.6);
  \draw[fill] (15,30) circle (0.6);
  \draw[fill] (0,35) circle (0.6);
  \draw[line width=.7pt] (0,5)--(-15,10)--(15,20)--(-15,30)--(0,35)--
                         (15,30)--(-15,20)--(15,10)--(0,5);
  \draw[line width=.7pt] (-15,10)--(-15,30);
  \draw[line width=.7pt] (15,10)--(15,30);
\end{tikzpicture}
}
    \end{picture}
  \caption{The Tardos poset}  
  \end{figure}  

Later, Z\'adori~\cite{Zadori}
studied series-parallel posets, a fairly large class of
posets that are very far from lattices.
For example, the Tardos poset is a series-parallel poset.
Z\'adori's results imply that
if $(A;\le)$ is a finite bounded series-parallel poset, then $\Pol(\le)$
is finitely generated if and only if the Tardos poset is
$\underline{\text{not}}$ a retract
of $(A;\le)$, and this is the case
if and only if $\Pol(\le)$ contains a $5$-ary near unanimity operation.
An immediate consequence of this result is that for every finite set $A$ of
size $|A|\ge8$ there exists a bounded partial order $\le$ such that the maximal
clone $\Pol(\le)$ is not finitely generated.

\subsection{Submaximal clones}
The results discussed in the last two paragraphs 
all point in the same direction:
by continuing Rosenberg's work on the maximal clones,
and using his ideas and techniques
to determine the maximal subclones of (at least)
the finitely generated maximal clones,
we can gain a better understanding of the 
lattice of clones on arbitrary finite sets, at least `near the top'.
A clone on a fixed finite set $A$ ($|A|\ge3$) 
is called a \emph{submaximal clone} on $A$ if it is a 
maximal subclone of a maximal clone on $A$.
The concept was introduced by Rosenberg in \cite{Ros-submax}, where he proved
a S{\l}upecki-type completeness theorem for the maximal clones $\Pol(B)$
where $B$ is a unary central relation.
The collection of the submaximal clones on $A$
may be thought of as the `next level' in the clone lattice on $A$
below the `top level' of the maximal clones; however, it should be noted that
not all submaximal clones on $A$ have depth~$2$ below $\OO_A$, because
there exist submaximal clones that are contained in two maximal clones so that
they are maximal in one of them and not maximal in the other.

Over the last five decades, a large body of work has been published about
submaximal clones on finite sets $A$ ($|A|\ge3$),
but the project of finding all of them is far from complete.
In the special case $|A|=3$ all submaximal clones on $A$ are known; their
description was completed by Lau~\cite{Lau-submax-3el}, using earlier
results by Machida~\cite{Machida-3el-mon-submax},
Demetrovics, Hann\'ak, Marchenkov~\cite{Demetrovics-Hannak-Marchenkov}, and
Salomaa~\cite{Salomaa}
on the submaximal clones on $3$-element sets 
that are maximal in maximal clones of types (1), (2), and (3), respectively. 
A complete list of all submaximal clones is not known for any finite set
of size $>3$, but there are a lot of maximal clones
$\MM$ on arbitrary finite sets $A$ ($|A|\ge3$)
for which all of its maximal subclones are known.
These include the following, where the numbers $(1)$--$(6)$  
indicate the type of the maximal clone $\MM$:
\begin{enumerate}
\item[(1)]
  $\MM=\Pol(\le)$ where $(A;\le)$ is a chain and $|A|\le5$;
  Larose~\cite{Larose}.
\item[(2)]
  $\MM$ is arbitrary of type (2); Rosenberg--Szendrei~\cite{Ros-Sz-perm-submax}.
\item[(3)]
  $\MM$ is arbitrary of type (3); Szendrei~\cite{Sz-affine-submax}.
\item[(4)]
  $\MM=\Pol(\rho)$ where $\rho$ is an equivalence relation on $A$ with
  exactly two blocks, one of which is a singleton;
  Lau \cite[Sec.~18.3]{Lau-book}.
\item[(5)$_1$]
  $\MM=\Pol(\{b\})$, $b\in A$; Lau~\cite{Lau-pol_0-submax}.
  \begin{enumerate}
  \item[]
    More generally, the list of maximal subclones is known also
    for all intersections $\bigcap_{b\in B}\Pol(\{b\})$
    ($\emptyset\subsetneq B\subseteq A$) of maximal clones of this form;
    Szendrei~\cite{Sz-cap_pol_a-submax},
      Lau~\cite{Lau-cap_pol_a-submax}, cf.\ also \cite[Ch.~16]{Lau-book}.
  \end{enumerate}  
  \item[(5)$_1$]
    $\MM=\Pol(B)$, $B\subsetneq A$, $|B|>1$; Lau~\cite[Ch.~17]{Lau-book}.
  \item[(6)$_{*}$]
    $\MM$ is S{\l}upecki's clone on $A$; Szendrei~\cite{Sz-slup-submax}.
  \begin{enumerate}
  \item[]
    In \cite{Sz-slup-submax},
    the list of maximal subclones is given also
    for every subclone of S{\l}upecki's clone on $A$,
    which contains all non-surjective operations on $A$.
  \end{enumerate}  
\end{enumerate}  

\subsection{Minimal varieties generated by finite algebras}
Rosenberg's Completeness Theorem has also been
instrumental in solving some problems that are apparently unrelated
to completeness. One of these is the problem 
to characterize the minimal
varieties that contain a nontrivial finite algebra.
Each such variety is 
of the form $\var{V}(\al{A})$
for a nontrivial finite algebra $\al{A}$,
and by choosing a generator $\al{A}$ of minimal cardinality, we get
that the variety has the form $\var{V}(\al{A})$ for a  
finite, strictly simple algebra $\al{A}$.
Thus, the problem may be rephrased
as stated in Problem~10 in \cite[p.~192]{Hobby-McKenzie}:
describe all finite, strictly simple algebras that generate minimal varieties.
Most strictly simple algebras $\al{A}$ that are `complete' in one of the
senses we discussed (e.g., primal, quasi-primal, or para-primal)
are known to generate minimal varieties;
the only exceptions are
the strictly simple affine algebras that have no trivial subalgebras
(see
\cite{Foster-boolean2,
  Pixley-quasiprimal, Clark-Krauss1, McKenzie-ms1976,
  McKenzie-paraprimal}).
At the other extreme, the $2$-element algebras $\al{A}=(A;\emptyset)$
with no operations (essentially, just $2$-element sets)
and the $2$-element algebras $\al{A}=(A;c_a)$ with a single constant operation
$c_a$ with value $a\in A$ (pointed sets),
are also strictly simple and generate minimal varieties, but  
their clones are too small to consider them `complete' 
is any reasonable sense.

The problem stated in the preceding paragraph
was solved by Kearnes and Szendrei
in~\cite{KSz-minvar-locfin} by the following theorem:

\smallskip
\noindent
{\it
A finite, strictly simple algebra $\al{A}$ generates a minimal
variety if and only if
\begin{enumerate}
\item[(i)]
  $\al{A}$ is either non-abelian or has a trivial subalgebra, and
\item[(ii)]
  for any (equivalently, for some) unary operation $e$
  in $\Clo(\al{A})$
  such that
  \begin{equation}\label{eq-min-idemp}
    \text{%
      $e=e^2$, $e$ is not constant, but $|e(A)|$ is as small as possible,
    }
  \end{equation}
  there exist operations $f_0,\dots,f_n$ in $\Clo^{(2)}(\al{A})$ and
  $g_0,\dots,g_n$, $h_0,\dots,h_n$
  in $\Clo^{(1)}(\al{A})$, for some $n\ge1$, such that $\al{A}$ satisfies
  the following identities:
  \begin{multline*}
    \qquad\qquad
    x = f_0 (x, eg_0(x)),\ \ 
    f_{i-1} (x, eh_{i-1} (x)) = f_i (x, eg_i (x))\ (1\le i\le n),\\ 
    f_n (x, eh_n (x)) = e(x).\ \ 
  \end{multline*}
\end{enumerate}}  
\noindent
It follows also from this characterization that in every minimal variety
$\var{V}(\al{A})$ where $\al{A}$ is a finite, strictly simple algebra,
the generating algebra
$\al{A}$ embeds in every nontrivial member of $\var{V}(\al{A})$ and hence 
$\al{A}$ is the only strictly simple algebra in $\var{V}(\al{A})$,
up to isomorphism.
In the proof, the operation $e\in\Clo^{(1)}(\al{A})$ is used to construct
from $\al{A}$ its \emph{induced algebra}
$e(\al{A})=\bigl(e(A);\{ef|_{e(A)}\colon f\in\Clo(\al{A})\}\bigr)$,
an algebra on the range
$e(A)$ of $e$ whose operations are the restrictions to $e(A)$ of
all operations in $\Clo(\al{A})$ that map into $e(A)$.
It is not hard
to check that $\al{M}:=e(\al{A})$ inherits many relevant properties
of $\al{A}$: for example, (a)~since $\al{A}$ is strictly simple,
so is $\al{M}$, (b)~$\al{M}$ is abelian if and only if $\al{A}$ is,
and (c)~$\al{M}$ has the same trivial subalgebras as $\al{A}$.
Furthermore, by the choice of $e$ described in \eqref{eq-min-idemp},
(d)~every unary operation in $\Clo(\al{M})$ is either a permutation
or a constant. A finite algebra $\al{M}$ with property~(d) is called
\emph{term minimal}.
The main idea of the proof of the characterization above
is to establish first that
for a strictly simple, term minimal algebra $\al{M}$, the variety
$\var{V}(\al{M})$ is minimal if and only if 
condition (i) holds for $\al{M}$ (in place of $\al{A}$), and then to
determine the exact relationship between $\al{A}$ and
its induced term minimal algebra $\al{M}=e(\al{A})$ that
allows the minimality of $\var{V}(\al{M})$ to be pulled back to
yield the minimality of $\var{V}(\al{A})$.
The first step here relies on a detailed analysis of all the possibilities
for the clones $\Clo(\al{M})$ of strictly simple, term minimal algebras
$\al{M}$. This analysis was completed
in \cite{Sz-termmin}, relying heavily on earlier results from
\cite{Post} for the case $|M|=2$ and from
\cite{Palfy, Sz-idemp, Sz-ss-surjective} for the case $|M|\ge3$.

A strengthening of Rosenberg's Completeness Theorem
from \cite{Sz-primal-alg-char-thm} plays a crucial role in this analysis
in the case where $\al{M}$ is a strictly simple, term minimal algebra
such that $\Clo^{(1)}(\al{M})$ is a transitive permutation group on $M$
($|M|\ge3$).
If, for the time being, we let
$\al{M}$ be an arbitrary finite, strictly simple algebra, then we have that
$\Clo(\al{M})\not\subseteq\Pol(\rho)$ holds for the Rosenberg relations
$\rho$ on $M$ of type (4) and (5)$_1$. Therefore, for such algebras $\al{M}$,
Rosenberg's Completeness Theorem can be restated as follows:
either $\al{M}$ is primal or $\Clo(\al{M})\subseteq\Pol(\rho)$ holds
for one of the Rosenberg relations $\rho$ on $M$ of type (1), (2), (3),
non-unary (5), or (6).
The strengthening of this theorem in \cite{Sz-primal-alg-char-thm}
replaces the condition ``$\Clo(\al{M})\subseteq\Pol(\rho)$ for $\rho$ of type
(2) or (3)'' by the following, more restrictive condition: either
the algebra
\begin{enumerate}
\item[(i)]
  $\al{M}$ is quasi-primal (but not primal) or affine, or
\item[(ii)]
  $\al{M}$ is isomorphic to an algebra $\al{M}'=(\mathbf{2}^r;F')$
  where $\mathbf{2}=\{0,1\}$, $r\ge2$, and $F'$ --- and hence
  also $\Clo(\al{M}')$ --- 
consists entirely of
`selector functions', that is, operations of the form
described in \eqref{eq-reg-max-clones}
where $g_1,\dots,g_r$ are unary operations on the set $\mathbf{2}$.
\end{enumerate}
Note that if $\Clo(\al{M})$ satisfies this last condition, then
$\Clo(\al{M})\subseteq\Pol(\rho)$ holds
for a Rosenberg relation $\rho$ on $M$ of type (3). 
In the special case when $\al{M}$ also has the property that
$\Clo^{(1)}(\al{M})$ is a transitive permutation group on $M$, then
it follows (see \cite{Sz-ss-surjective}) that $\al{M}$ is quasi-primal or
affine or essentially unary.

\subsection{The Dichotomy Theorem for Constraint Satisfaction
Problems}\label{subsec-csp}
To conclude this section,
we will discuss a problem from theoretical computer science
where Rosenberg's Completeness Theorem was key to a solution. The
problem is the Dichotomy Conjecture for nonuniform, finite domain constraint
satisfaction problems (briefly: CSPs), and the solution we will discuss is 
Zhuk's proof of the conjecture in
\cite{Zhuk-dich-conj_conf, Zhuk-dich-conj_full}.\footnote{%
Simultaneously and independently, Bulatov also proved the Dichotomy
Conjecture, see \cite{Bulatov}.
His approach is quite different from Zhuk's, and does not use Rosenberg's
Completeness Theorem.}
Nonuniform, finite domain CSPs form a large class of combinatorial
decision problems, with a vast literature in computer science,
which includes many familiar problems, such as Boolean satisfiability,
$k$-colorability of graphs, scheduling problems, and solving
systems of linear equations over finite fields.
They are parameterized by finite relational structures
$\str{A}=(A;\rho_1,\dots,\rho_m)$,
where the relations $\rho_i$ are nonempty,
and the size of the domain $A$ as well as the number of
relations $\rho_i$ are finite.
The decision problem 
$\CSP(\str{A})$ is the following:
given any finite relational structure
$\str{X}=(X;\rho_1^{\str{X}},\dots,\rho_m^{\str{X}})$ as an input,
where for each $i$ the relation
$\rho_i^{\str{X}}$ has the same arity $k_i$ as the corresponding
relation $\rho_i$ of
$\str{A}$, the task is to
determine whether or not there exists a homomorphism
$v\colon\str{X}\to\str{A}$.
Here, a homomorphism $v\colon\str{X}\to\str{A}$ is
a map $v\colon X\to A$
such that whenever a tuple $(x_1,\dots,x_{k_i})\in X^{k_i}$ is in
$\rho_i^{\str{X}}$ then the image tuple
$(v(x_1),\dots,v(x_{k_i}))\in A^{k_i}$ is in $\rho_i$.
For applications, one should think of the elements of $X$ as `variables',
to which we want to assign `values' from the set $A$ via 
$v\colon X\to A$; however, $v$ has to satisfy some `constraints':
each pair $\bigl((x_1,\dots,x_{k_i}),\rho_i\bigr)$ with
$(x_1,\dots,x_{k_i})\in \rho_i^{\str{X}}$ is a constraint which requires
that the values assigned to these variables have to satisfy
the condition $(v(x_1),\dots,v(x_{k_i}))\in \rho_i$.
The Dichotomy Conjecture, due to
Feder and Vardi~\cite{Feder-Vardi1993, Feder-Vardi1998},
states that
for every finite relational structure $\str{A}$, the problem
$\CSP(\str{A})$ is either in P (i.e., solvable by a deterministic
algorithm in polynomial time) or is NP-complete.\footnote{%
It is easy to see that each such problem $\CSP(\str{A})$ is in NP, the
point of the conjecture is that $\CSP(\str{A})$ cannot be of intermediate
complexity.}

Let $\str{A}=(A;\rho_1,\dots,\rho_m)$ be a
finite relational structure, and
define the \emph{polymorphism clone} $\Pol(\str{A})$
of $\str{A}$ to be the clone $\bigcap_{i=1}^m\Pol(\rho_i)$
consisting of all operations on $A$ that preserve every relation of
$\str{A}$.
It was noticed by Jeavons, Cohen, and Gyssen (see \cite{Jeavons-C-G, Jeavons})
that if $\str{A}=(A;\rho_1,\dots,\rho_m)$ and
$\str{A'}=(A;\sigma_1,\dots,\sigma_n)$ are finite relational
structures on the same set $A$
such that, under the Galois correspondence between operations and
relations on $A$ (see \cite{BKKR, Gei}),
the Galois-closed set of relations generated by
the relations of $\str{A}$ contains the relations of $\str{A}'$, then
the problem $\CSP(\str{A}')$ is log-space (hence polynomial-time)
reducible to $\CSP(\str{A})$. Thus, up to log-space equivalence,
the computational complexity of $\CSP(\str{A})$ depends only on the
clone $\Pol(\str{A})$.
Bulatov, Jeavons, and Krokhin showed (see
\cite{Bulatov-Jeavons-Krokhin2000, Bulatov-Jeavons,
  Bulatov-Jeavons-Krokhin2005})
that if $e$ is a unary operation in $\Pol(\str{A})$ satisfying
the conditions in \eqref{eq-min-idemp}, except that we allow $e$ to be constant,
then the relational structure $e(\str{A})=(e(A);e(\rho_1),\dots,e(\rho_m))$
obtained from $\str{A}$ by restricting the base set as well as all the
relations to the range $e(A)$ of $e$ has the following properties:
\begin{enumerate}
\item[(a)]
  $e(\str{A})$ is a \emph{core structure}, meaning that every homomorphism
  $e(\str{A})\to e(\str{A})$ $\bigl($i.e., every unary operation in
  $\Pol(e(\str{A}))\bigr)$ is a permutation, and
\item[(b)]
  $\CSP(\str{A})$ and $\CSP(e(\str{A}))$ are essentially the same
  problem, that is, they have the same `YES/NO' answer for every input, because
  $\str{A}\stackrel{e}{\to}e(\str{A})$
  and the inclusion map $e(\str{A})\stackrel{\textrm{incl}}{\to}\str{A}$
  are both homomorphisms.
\end{enumerate}  
Therefore, to prove or disprove the Dichotomy Conjecture,
no generality is lost by restricting to problems $\CSP(\str{A})$ where
$\str{A}$ is a core structure.
We can also assume $|A|\ge2$, because
in the case $A=\{a\}$ each relation $\rho_i$ of $\str{A}$ is a singleton
$\{(a,\dots,a)\}$, and the answer to the decision problem $\CSP(\str{A})$
is `YES' for every input
$\str{X}$, as witnessed by the constant homomorphism
$v\colon\str{X}\to\str{A}$ with value $a$.
Bulatov, Jeavons, and Krokhin also showed that
\begin{enumerate}
\item[(c)]
if $\str{A}$ is a core structure with $|A|\ge2$, then
the problem $\CSP(\str{A})$ is log-space equivalent to the problem
$\CSP(\str{A}^*)$ where $\str{A}^*$ is the relational structure
obtained from $\str{A}$ via expanding it by all singleton relations 
$\{a\}$ ($a\in A$).
\end{enumerate}  
The effect of the construction $\str{A}\rightsquigarrow\str{A}^*$
on the polymorphism clones is that
$\Pol(\str{A}^*)$ is the subclone of $\Pol(\str{A})$ consisting of all
\emph{idempotent} operations $f\in\Pol(\str{A})$, that is,
all operations $f\in\Pol(\str{A})$ which satisfy the identity
$f(x,\dots,x)=x$.

By analyzing the cases for which the complexity of $\CSP(\str{A})$ was
known at the time,
Bulatov, Jeavons, and Krokhin strengthened the Dichotomy Conjecture of
Feder and Vardi 
to a conjecture that also predicts the dividing line between
the problems $\CSP(\str{A})$ that are in P and those that are NP-complete
(assuming P${}\not={}$NP).
This conjecture became known as the Algebraic Dichotomy Conjecture. It has
a number of equivalent formulations, one of which is the
following:
\begin{equation}\label{eq-AlgDichConj}
\begin{matrix}
  \text{\it If $\str{A}$ is a finite relational structure with at least
    two elements and}\hfill\\
  \text{\it $\str{A}$ is a core, then}\hfill\\
  \text{\it \phantom{m}$(\lozenge)$\phantom{i}
  $\CSP(\str{A})$ is in P if $\Pol(\str{A})$
    contains a Taylor operation, and}\hfill\\
  \text{\it \phantom{m}$(\blacklozenge)$\phantom{i}
  $\CSP(\str{A})$ is NP-complete otherwise.}\hfill
\end{matrix}
\end{equation}

A \emph{Taylor operation} (named after Taylor~\cite{Taylor})
on a set $A$ is an idempotent operation $f$
on $A$, which satisfies strong enough identities that prevent $f$ from
being a projection operation.
\comment{%
say $n$-ary ($n\ge2$),
such that for every $\ell$ ($1\le \ell\le n$), $f$ satisfies an identity
of the form
\begin{equation*}
  f(\dots,\underbrace{x}_\ell,\dots)=
  f(\dots,\underbrace{y}_\ell,\dots)
\end{equation*}  
where $x$ and $y$ are distinct variables in the $\ell$-th argument of $f$,
and all other variables are $x$ or $y$.}%
For example, Mal'tsev operations, majority operations, and near
unanimity operations are Taylor operations by their defining identities
in \eqref{eq-maltsev}, \eqref{eq-majority}, \eqref{eq-nu},
and so are all binary operations
$f$ satisfying the identity $f(x,y)=f(y,x)$.
The easy part, $(\blacklozenge)$, of the conjecture was established in
\cite{Bulatov-Jeavons-Krokhin2000, Bulatov-Jeavons,
  Bulatov-Jeavons-Krokhin2005}. 
To prove $(\lozenge)$, it is useful to strengthen the condition
``$\Pol(\str{A})$ contains a Taylor operation'' in $(\lozenge)$ to 
``$\Pol(\str{A})$ contains a specific type of Taylor operation'' where
the identities ensuring the Taylor property are more manageable.
Several such strengthenings emerged over the years. The one
Zhuk uses relies on the following theorem
due to Mar\'oti and McKenzie~\cite{Maroti-McKenzie}:

\smallskip
\noindent
{\it
If a clone on a finite set contains a Taylor operation,
then it also contains
a \emph{special weak near unanimity operation}
(briefly: \emph{special WNU operation})
$u$, which is an idempotent operation $u$
that satisfies {\rm(i)}~the identities obtained from the near unanimity
identities in \eqref{eq-nu} by deleting ``$=x$'', and {\rm(ii)}~the identity
\begin{equation*}
  u(x,\dots,x,u(x,\dots,x,y))=u(x,\dots,x,y).
\end{equation*} }

These preparations show that to prove $(\lozenge)$, it suffices to produce
a polynomial time algorithm for solving $\CSP(\str{A})$ for any fixed
finite structure $\str{A}$ such that $\Pol(\str{A})$ contains a
special WNU operation $u$, and the relations of $\str{A}$ are
all
relations on its base set $A$ that are preserved
by $u$ and have arity $\le k$ for some fixed integer $k\ge2$.
Zhuk proved the Algebraic Dichotomy Conjecture in 
\cite{Zhuk-dich-conj_conf, Zhuk-dich-conj_full} by presenting such
an algorithm.
The main algorithm is recursive, with recursive calls applied to
instances that are in some well-defined sense `smaller' or `simpler' than the
instances at hand.
The main algorithm starts with some combinatorial preprocessing
of the input structure $\str{X}$, which has one of the following outcomes:
(a)~it is decided that no homomorphism $\str{X}\to\str{A}$ exists, or
(b)~a `smaller' or
`simpler', but equivalent input structure $\str{X}'$ is produced
(where `equivalent' means:
a homomorphism $\str{X}\to\str{A}$ exists if and only if 
a homomorphism $\str{X}'\to\str{A}$ exists),
which allows the main algorithm to be called recursively for
$\str{X}'$ in place of $\str{X}$, or
(c)~we know that
$\str{X}$ has strong consistency and irreducibility properties. 
From this point on the main algorithm continues by 
assuming that $\str{X}$ is as in case~(c), and performs further reductions
that form the main body of the algorithm.
These reductions are algebraic in nature, and are
based on the following consequence of
Rosenberg's Completeness Theorem:

\smallskip
\noindent
{\it
If $\al{D} = (D; w)$ is a finite algebra where $w$ is an $m$-ary
special WNU operation, then one of the following
conditions holds:
\begin{enumerate}
\item[(i)]
  there exists a proper subalgebra $\al{B}$ of $\al{D}$
  and a binary operation
  $t\in\Clo(\al{D})$ such that
  $t(D,B)\cup t(B,D)\subseteq B$;
\item[(ii)]
  there exists a proper subalgebra $\al{C}$ of $\al{D}$, another
  algebra $\al{H}=(H;w_H)$ with an $m$-ary
  special WNU operation $w_H$, and a subdirect
  subalgebra $\al{R}$ of $\al{D}\times\al{H}$ such that condition~(i)
  fails for $\al{H}$ (in place of $\al{D}$), and
  \begin{equation*}
    C=\{c\in D:\{c\}\times H\subseteq R\};
  \end{equation*}
\item[(iii)]  
  there exists a maximal congruence $\sigma$ on $\al{D}$ such that 
  $\al{D}/\sigma$ is functionally complete;
\item[(iv)]
  there exists a maximal congruence $\sigma$ on $\al{D}$ such that
  $\al{D}/\sigma$ is isomorphic to $(\al{Z}_p; x_1+\dots+ x_m )$
  for some prime divisor $p$ of $m-1$.
\end{enumerate}}

\smallskip
\noindent
Applying this theorem to the subalgebras $\al{D}_x$ of $(A;u)$
associated to each element $x$ of the current input
$\str{X}$ (where the value $v(x)$
of a potential homomorphism $v\colon\str{X}\to\str{A}$ may lie),
Zhuk proves that in each one of the cases (i)--(iii),
a `smaller', but equivalent input structure $\str{X}''$ can be produced
from $\str{X}$, which
again allows the main algorithm to be called recursively for
$\str{X}''$ in place of $\str{X}$.
Handling the remaining case, when each subalgebra $\al{D}_x$ ($x\in X$)
of $(A;u)$ is as in case (iv) for some prime $p_x$,
is the most difficult part of the main algorithm, and requires developing
and using new algebraic tools, some akin to results in commutator theory.

\section{Maximal clones in other contexts}
Rosenberg pioneered the work on maximal clones in several other contexts.

\subsection{Maximal clones on infinite sets}
Rosenberg proved in \cite{Ros-number-maxclones-infsets} that there are
$2^{2^{|A|}}$ maximal clones on any infinite set $A$. This is equal to
the number of all clones on $A$,
which is in sharp contrast to what is true on finite base sets.
The proof follows the basic idea of Gavrilov's paper \cite{Gavrilov}
on the number of maximal clones on countably infinite sets $A$, and
demonstrates the existence of $2^{2^{|A|}}$ maximal clones on arbitrary infinite
sets $A$, without exhibiting any maximal clone. Later,
Goldstern and Shelah~\cite{Goldstern-Shelah-constr-maxclones-infsets}
found a transparent explicit construction, which produces
different maximal clones from different ultrafilters on infinite sets $A$,
yielding a simpler proof of Rosenberg's result.   

In \cite{Ros-number-maxclones-infsets}, Rosenberg also
brought up the following question about clones on infinite sets $A$:  
\begin{align}
\text{Ques}&\text{tion: Is it true that}\notag\\
&\text{every proper subclone of $\OO_A$ is contained in a maximal clone on $A$?}
\label{eq-dually-atomic}
\end{align}
We saw that for finite sets $A$, the answer is `yes'.
In contrast, 
for infinite sets $A$, this turned out to be a surprisingly
difficult problem: to this day, no infinite set $A$ is known for which
a final answer `yes' or `no' has been found.
Note that by saying that the answer
is `yes' we mean that the statement in \eqref{eq-dually-atomic} can be proved
in ZFC\footnote{%
ZFC is the abbreviation for ``Zermelo--Fraenkel axioms, together with
the Axiom of Choice''.},
the usual axiom system for set theory, while by saying the answer
is `no' we mean that
the negation of the statement in \eqref{eq-dually-atomic} can be proved
in ZFC.
It is easy to see that the answer to this question depends only on the
cardinality of the set $A$,
since the clone lattices on sets of the same cardinality are isomorphic.
The strongest result concerning this question was obtained by
Goldstern and Shelah in
\cite{Goldstern-Shelah-clonelattice-not-duallyatomic_countable-base} and
\cite{Goldstern-Shelah-clonelattice-not-duallyatomic_some-uncountable-bases},
where the authors proved the following theorem
first for the case when the cardinality, $\kappa$, of $A$ is $\aleph_0$
(i.e., $A$ is countably infinite) and more recently
for the case when $\kappa$ is an uncountable regular cardinal:

\smallskip
\noindent
{\it
The negation of the statement in~\eqref{eq-dually-atomic}
for sets $A$ of regular cardinality $\kappa$
is a consequence of
the axiom system where ZFC is expanded by the assumption $\kappa^+=2^{\kappa}$.
}

\smallskip
\noindent
The assumption $\kappa^+=2^{\kappa}$ here means that
the Generalized Continuum Hypothesis is assumed to hold at
the cardinal $\kappa$.
The proof of the theorem
proceeds by constructing a clone $\CC$ on $A$, and a subfamily
$\mathsf{F}$
of the interval $[\CC,\OO_A)$ (containing $\CC$ but not
containing $\OO_A$) such that (a)~$\mathsf{F}$ has no maximal elements, but
(b)~$\mathsf{F}$ is \emph{cofinal} in $[\CC,\OO_A)$ (with respect to inclusion),
that is, every clone in $[\CC,\OO_A)$ is a subclone of a member of
$\mathsf{F}$. 
Since ZFC\footnote{%
Assuming, as usual, that ZFC is consistent.}
has models where $\kappa^+=2^\kappa$ holds for every
infinite cardinal $\kappa$, this theorem shows that for infinite sets $A$ of
regular cardinality,
\begin{itemize}
\item
  the negation of \eqref{eq-dually-atomic} holds in some models of set theory,
  and hence
\item
  the statement in \eqref{eq-dually-atomic} is not provable in ZFC.
\end{itemize}  
For infinite sets $A$ of singular cardinality, the problem
is completely open.\footnote{%
For more details about clones on infinite sets, the reader is referred to
the survey paper~\cite{Goldstern-Pinsker-survey}.} 

\subsection{Completeness for locally closed clones on infinite sets}
If $\CC$ is a clone on a set $A$, an operation $f$ on $A$ is said to be
\emph{locally in $\CC$} if on any finite subset of its domain, $f$ agrees 
with, i.e., $f$ can be \emph{interpolated by}, some member of $\CC$
(of the same arity as $f$). Furthermore, 
$\CC$ is called \emph{locally closed} if
it contains every operation on $A$ that is locally in $\CC$.
Of course, if $A$ is finite, then every clone
on $A$ is locally closed, so these notions are interesting only if the
base set $A$ is infinite.
It is easy to see that a clone $\CC$ on $A$ is locally closed if and only if
for every $n\ge1$, the $n$-ary part $\CC^{(n)}$ of $\CC$ is closed in the
product topology of $A^{A^n}$ where $A$ is equipped with the discrete topology.
For every clone $\CC$ on $A$ there exists a least locally closed clone
containing $\CC$, which is called the \emph{local closure of $\CC$}, and
consists exactly of those operations on $A$ that are locally in $\CC$.
The local closure of the clone $\Clo(\al{A})$ of an algebra $\al{A}=(A;F)$,
that is, the least locally closed clone on $A$ containing $F$,
is called the \emph{locally closed clone of $\al{A}$}. 
For a classical example, consider the clone
$\mathcal{P}_{\mathbb{R}}$
of all polynomial functions
with real coefficients in finitely many variables on the set $\mathbb{R}$
of real numbers.
By Lagrange's Interpolation Theorem, the local closure of this clone
is the clone $\OO_{\mathbb{R}}$ of all functions
$\mathbb{R}^n\to\mathbb{R}$ ($n\ge1$). 
This example also shows that on an infinite base set,
the local closure of a clone may be much bigger
than the clone itself.
For example, by looking at the the cardinalities of 
the clone $\mathcal{P}_{\mathbb{R}}$
of polynomial functions on $\mathbb{R}$
and its local closure $\OO_{\mathbb{R}}$, we see that 
$\mathcal{P}_{\mathbb{R}}$
has cardinality
$|\mathbb{R}|=2^{\aleph_0}$, while $\OO_A$ 
has cardinality
$2^{|\mathbb{R}|}=2^{2^{\aleph_0}}$.

Foster noticed in \cite{Foster-loc-func-completeness} --- much before the
concept of locally closed clones was isolated and studied --- that
many of the results he proved earlier about the structure of algebras
in the variety $\var{V}(\al{A})$ generated by a primal or functionally
complete algebra $\al{A}$ carry over to infinite algebras that are
primal or functionally complete `locally'.\footnote{He used the phrase
primal or functionally complete `in the small'.}
In \cite{Pixley-discriminator}, Pixley considered locally
quasi-primal algebras along with quasi-primal algebras,
and carried over his characterizations of
quasi-primal algebras (see Subsection~\ref{subsec-quasiprimal})
to locally quasi-primal algebras.
The Baker--Pixley Theorem was also proved in
\cite{Baker-Pixley} for infinite algebras $\al{A}$ as well.
In all these extensions from the finite to the infinite case,
the only change needed in both the
definitions and the theorems was
that the clone of $\al{A}$ had to be replaced
by the locally closed clone of $\al{A}$.

Why do these completeness-type results for finite algebras
carry over so smoothly to infinite algebras in a `local version'?
The reason was revealed
in the late 1970's by studying the Galois connection between
(finitary) operations on infinite sets $A$ and finitary relations 
on $A$. Extending the results of \cite{BKKR}, Romov announced in
\cite{Romov-finfin-galoisconn-infsets}
that in this Galois connection the Galois-closed sets of operations
are exactly the locally closed clones.\footnote{For more about
this Galois connection, including the description of the Galois-closed sets
of relations, and its applications to concrete representations of
related structures of algebraic structures, see
Szab\'o~\cite{Szabo-finfin-galoisconn-infsets} and   
P\"oschel~\cite{Poschel-finfin-galoisconn-infsets1,
  Poschel-finfin-galoisconn-infsets2}.}
Hence, every locally closed clone on any set $A$ is
of the form $\Pol(\mathsf{R}):=\bigcap_{\rho\in\mathsf{R}}\Pol(\rho)$
for some set $\mathsf{R}$ of finitary relations on $A$.

In particular, every clone that is maximal among the locally closed clones on
any set $A$ has the form
$\Pol(\rho)$ for some finitary relation $\rho$ on $A$.
Rosenberg and Schweigert were the first to study these clones,
which they called \emph{locally maximal clones},
in \cite{Ros-Schw-loc-max}.
They proved that $\Pol(\rho)$ is locally
maximal on $A$ for each one of the following relations, which are
recognizably very close to
some of the relations familiar from Rosenberg's Completeness Theorem:
\begin{enumerate}
\item[(1)$'$]
  $\rho$ is a locally bounded partial order on $A$;
\item[(2)]
  $\rho$ is (the graph of) a fixed point free permutation on $A$
  such that all cycles of $\rho$ have the same prime length;
\item[(3)$'$]
  $\rho$ is the affine relation $\alpha^+$ obtained from an
  abelian group $(A;+)$, as described in \eqref{eq-affine-rel}, 
  provided $(A;+)$ is either an elementary abelian $p$-group
  for some prime $p$ or torsion-free and divisible;
\item[(4)]
  $\rho$ is a nontrivial equivalence relation on $A$;
\item[(5)$'$]
  $\rho$ is a locally central relation on $A$.
\end{enumerate}  
Let $\mathsf{R}_1$ denote the set of all these relations on $A$.

It was also proved in \cite{Ros-Schw-loc-max} that the analog of
Statement~\eqref{eq-dually-atomic}
is false for locally closed clones on infinite sets $A$.
For example, if $\rho$ is (the graph of) a fixed point free permutation
on $A$ such that all cycles of $\rho$ are infinite, then
$\Pol(\rho)$ is not contained in any locally maximal clone on $A$.
Therefore, to obtain a completeness criterion for $\OO_A$, as a locally
closed clone on an infinite set $A$, it is not enough to find all locally
maximal clones on $A$.
In \cite{Ros-Szabo-loc-max}, Rosenberg and Szab\'o isolated several other types
of finitary relations --- one of the types being the graphs of fixed point
free permutations where all cycles are infinite ---,
and proved that for the set $\mathsf{R}_2$
of all these relations, the family
$\{\Pol(\rho):\rho\in\mathsf{R}_1\cup\mathsf{R}_2\}$ of locally closed clones
on $A$ is cofinal in the set of all proper, locally closed subclones of
$\OO_A$ (ordered by inclusion). 
This implies the following completeness theorem:

\smallskip
\noindent
{\it A set $F$ of functions on $A$
is \emph{locally complete}, that is, the local closure of the clone generated by $F$ is
$\OO_A$, if and only if $F\not\subseteq\Pol(\rho)$ holds for all relations
$\rho\in\mathsf{R}_1\cup\mathsf{R}_2$.} 

\subsection{Maximal clones of partial
  operations}\label{subsec-max-partial-clones}
Clones can be formed from partial operations just as well as from
total operations: a set $\UU$ of partial operations on a fixed base set
$A$ is a \emph{clone of partial operations} on $A$, or briefly a
\emph{partial clone} on $A$, if $\UU$ is closed under composition and contains
the (total) projection operations. A \emph{strong partial clone} on $A$ is a
partial clone on $A$ that is closed under
the process of restricting its members
to arbitrary subsets of their domain. It follows that
every partial clone $\UU$ is contained in a least strong partial
clone, its \emph{strong closure}, which consists of all restrictions of
the members of $\UU$
to arbitrary subsets of their domain. Hence, a maximal partial clone $\MM$
on $A$ is a strong partial clone, unless its strong closure is
the clone $p\OO_A$ of all partial operations on $A$. 
Haddad, Rosenberg, and Schweigert proved in
\cite{Haddad-Ros-Schw-max-partial-clones_3}  
that on any base set $A$, there is a unique maximal partial clone
whose strong closure is $p\OO_A$, namely the partial clone
$\MM_0:=\OO_A\cup\{\emptyset\}$ consisting of all totally defined operations
on $A$ and the partial operation with empty domain. As a consequence, they
obtained a S{\l}upecki-type completeness theorem for partial operations
on finite sets $A$.

A Rosenberg-type completeness theorem for partial operations on finite sets
was obtained by
Haddad and Rosenberg. The results were first announced in
\cite{Haddad-Ros-max-partial-clones_0},  
and then published with detailed proofs in
\cite{Haddad-max-partial-clones_1}, \cite{Haddad-Ros-max-partial-clones_2},
and \cite{Haddad-Ros-max-partial-clones_4}.
For the special cases when $|A|=2$ or $|A|=3$ these
results were proved earlier by
Freivald~\cite{Freivald-max-partial-clones-on-2}, and
by Lau~\cite{Lau-max-partial-clones-on-3} and
Romov~\cite{Romov-max-partial-clones-on-3}\footnote{Three maximal partial
clones are missing from the list of maximal partial clones in this paper.},
respectively. 
Haddad and Rosenberg
showed that if $A$ is a finite set ($|A|\ge2$), then
Statement~\eqref{eq-dually-atomic}
remains true for the maximal subclones of the clone $p\OO_A$ of all partial
operations on $A$,
and they determined all maximal partial clones on $A$.
Since $\MM_0$ is the only maximal partial clone that is not a strong
partial clone, the main task was to describe the maximal partial clones
that are strong partial clones. For this, they used Romov's result 
in \cite{Romov-relational-descr-for-partial-clones} that
strong partial clones can be described by relations, that is,
every strong partial clone on a finite set $A$ is of the form
$\pPol(\mathsf{R}):=\bigcap_{\rho\in\mathsf{R}}\pPol(\rho)$ for some set
$\mathsf{R}$ of finitary relations on $A$, where
$\pPol(\rho)$ is the partial clone consisting of all partial operations
that preserve the relation $\rho$.
In particular, maximal partial clones are
of the form $\pPol(\rho)$ for a single relation $\rho$.
Haddad and Rosenberg described four types of relations $\rho$ on finite sets
$A$, all of arity $\le\max(4,|A|)$, such that the strong maximal partial clones
on $A$ are exactly the partial clones $\pPol(\rho)$ where $\rho$ is a relation
of one of the four types. The relations that occur in this list are quite
different
(but not disjoint) 
from the list of relations in Rosenberg's Completeness Theorem.
For example, for $2$-element sets $A$ there are $5$ maximal clones in
$\OO_A$ and $7$ maximal partial clones in $p\OO_A$
(see \cite{Post, Freivald-max-partial-clones-on-2}),
while the corresponding numbers for $3$-element sets are
$18$ and $58$,
respectively
(see \cite{Yab, Lau-max-partial-clones-on-3}).

\section{Rosenberg's Theorem on minimal clones}

The aim of this section is to present
Rosenberg's theorem on minimal clones, 
which was announced at a
conference in 1983, and appeared in \cite{Ros-minclones} in 1986.
While Rosenberg's theorem on maximal clones on finite sets
solved a central problem of its time,
his theorem on minimal clones opened up 
decades of research that is still ongoing.

\subsection{Background}\label{subsec-minclone-background}
When we say that a clone is `minimal', we mean that the clone is minimal
for the property of being nontrivial. More formally,
a clone $\CC$ on a set $A$ is called \emph{minimal} if it has exactly
one proper subclone: the clone $\II_A$
of projections on $A$.
Equivalently, $\CC$ is minimal if and only if
it contains a non-projection and
for every operation
$f\in\CC$ that is not a projection, $\CC$ is generated by
$f$.
Saying that $\CC$ is generated by $f$ is equivalent to saying that  
$\CC$ coincides with the clone of the algebra $\al{A}_f:=(A;f)$.
From now on, we will often use the notation $\al{A}_f$ for the
algebra $(A;f)$ where $f$ is an arbitrary operation on $A$.

Since a minimal clone is required to contain a non-projection,
a minimal clone on $A$ does not exist if $|A|<2$.
Therefore we will assume from now on that $|A|\ge2$. However,
unless explicitly stated otherwise, we won't assume that $A$ is finite.

A systematic study of minimal clones on finite sets  
was initiated by P\"oschel and Kalu\v{z}nin
in their book \cite{PK}.
In Section~4.4 of the book the authors
discussed some of the fundamental facts
about minimal clones on finite sets, and in Problem~12 on page 120, they
posed the problem of classifying all minimal clones on finite sets.
To summarize some basic facts about minimal clones here, 
we will start with a list of examples of minimal clones.
In the examples, $\al{Z}_n$ denotes the set (or the ring)
of integers modulo $n$. We continue to use the notation
$\mathbf{2}=\{0,1\}$ for the common base set of all
clones of Boolean functions, but
when convenient, we may also consider clones on the set $\al{Z}_2$
as clones of Boolean functions.
For most examples of minimal clones below we cite the paper or book where
the minimality of the given clone was proved.

\smallskip
\noindent
$\underline{\text{Examples of minimal clones}}$:
\begin{enumerate}
\item[(i)]
  The three clones generated by the unary Boolean functions $c_0$, $c_1$
  (constants), and $\neg$ (negation), respectively, on $\mathbf{2}$.
\item[(ii)]
  The clone generated by the ternary addition $x_1-x_2+x_3$
  of the additive group $(\al{Z}_p;+)$
  (P{\l}onka~\cite{Plonka1}).
\item[(iii)]
  The clone generated by the binary operation $(1-p)x_1+px_2$ in the clone
  of the group $(\al{Z}_{p^2};+)$
  (P{\l}onka~\cite{Plonka1}).
\item[(iv)] 
  The clone of a \emph{rectangular band}\footnote{%
  A \emph{rectangular band} is a semigroup $(S;\cdot)$
  satisfying the identities $x^2=x$ and $xyx=x$. $(S;\cdot)$ is a
  \emph{left} [or \emph{right}] \emph{zero semigroup} if it satisfies
  the identity $xy=z$ [or $xy=y$, respectively], i.e.,
  its operation $\cdot$ is projection onto the first [second] variable.} 
  that is not a left or right zero semigroup.
  For example, $(\al{Z}_6;3x_1+4x_2)$ is a rectangular band.
\item[(v)]
  The two clones generated by the Boolean functions
  $\wedge$ (conjunction) and $\vee$ (disjunction), respectively,
  on $\mathbf{2}$.
\item[(vi)]
  The clone generated by the median operation
  $(x_1\wedge x_2)\vee(x_1\wedge x_3)\vee(x_2\wedge x_2)$ of a lattice
  (P\"oschel--Kalu\v{z}nin~\cite{PK}).
\item[(vii)]
  The clone generated by the ternary dual discriminator
  operation\footnote{The dual discriminator operation was introduced by Fried
  and Pixley in \cite{Fried-Pixley}.}
  $d$ defined by
  $d(a,b,c) = a$ if $a = b$ and $d(a,b,c)= c$ otherwise ($a,b,c\in A$)
  on any set $A$ of size $\ge2$
  (Cs\'ak\'any--Gavalcov\'a~\cite{Csakany-Gavalcova}). 
\item[(viii)]
  The clone generated by the following $k$-ary operation $\ell_k$
  on a $k$-element set $A$ ($k\ge3$):
  $\ell_k(a_1,\dots,a_k)=a_k$ if $\{a_1,\dots,a_k\}=A$ and
  $\ell_k(a_1,\dots,a_k)=a_1$ otherwise ($a_1,\dots,a_k\in A$)
  (Cs\'ak\'any--Gavalcov\'a~\cite{Csakany-Gavalcova}). 
\end{enumerate}    
By inspecting Post's description of all clones of Boolean functions
in \cite{Post}, one can see that the list above contains all the
$7$ minimal clones of
Boolean functions, namely: the clones in Examples~(i), (v), (ii) for $p=2$,
and (vi) for the $2$-element lattice $(\mathbf{2};\wedge,\vee)$ (which
coincides with the clone (vii) for the $2$-element set $A=\mathbf{2}$).

Since every minimal clone $\CC$ on a set $A$
is generated by a single operation, in most cases
it is convenient to fix a generating operation $f$ for $\CC$, and study
when such a clone $\CC=\Clo(\al{A}_f)$
is minimal.
For this, it is useful to
choose $f\in\CC$ so that $f$ has minimal arity among the
non-projections in $\CC$ (see e.g., \cite{PK}, \cite{Csakany-3el-min2}, and
\cite{Ros-minclones}).
This choice has the effect that
\begin{equation}\label{eq-identifying-vars}
\text{upon identifying any two variables
  in $f$ we get a projection.}
\end{equation}
In particular, if the arity of $f$ is $n\ge2$, then $f$ is idempotent.
In \cite{Swierczkowski}, \'Swierczkowski observed 
the following fact\footnote{The statement in \cite{Swierczkowski}
uses a different terminology.}
(cf.~also \cite[Thm.~4.4.6]{PK}):

\smallskip
\noindent
{\it
If an operation $f$ of arity $n\ge4$
has property \eqref{eq-identifying-vars}, then for some 
$i$ $(1\le i\le n)$, $f$ satisfies all identities
\begin{equation}\label{eq-semiproj}
  f(y_1,\dots,y_n)=y_i\quad
  \text{where two of the variables $y_1,\dots,y_n$ are equal}.   
\end{equation} } 

\smallskip
\noindent
These identities express that
upon identifying any two variables
in $f$, $f$ always turns into a
projection onto the same variable $x_i$ ($1\le i\le n$) of $f$.
An $n$-ary operation satisfying the identities in~\eqref{eq-semiproj}
is called an \emph{$n$-ary semiprojection onto the $i$-th variable}
if $n\ge3$ and $f$ is not a projection.

If $f$ has property \eqref{eq-identifying-vars} and
arity $n=3$, then there are four possibilities for $f$:
(a)~$f$ is a semiprojection onto one of its variables, or
(b)~$f$ is a majority
operation, that is, the identities in \eqref{eq-majority} hold for $m=f$, or 
(c)~$f$ is a \emph{minority operation}, that is, it satisfies the identities
\begin{equation}\label{minority}
f(x,y,y)=f(y,x,y)=f(y,y,x)=x,
\end{equation}  
or
(d)~the identities in \eqref{eq-arithm} hold either for $q=f$ or for an
operation $q$ obtained from $f$ by permuting variables.
In case~(d), $\Clo(\al{A}_f)$ also contains a majority operation,
namely the operation $q(x_1,q(x_1,x_2,x_3),x_3)$.

Thus, every clone on $A$ other than the clone $\II_A$
of projections, contains
one of the following five types of operations $f$:\footnote{%
The numbering of these types is not
as consistent in the literature as the numbering (1) through (6) of the types
of maximal clones. 
We chose a numbering that follows closely
the order in which these cases
are listed in \cite[Thm.~2.9]{Ros-minclones}.}
\begin{enumerate}
\item[(I)]
  $f$ is a unary operation, $f\notin\II_A$;
\item[(II)]
  $f$ is a binary idempotent operation, $f\notin\II_A$;
\item[(III)]
  $f$ is a semiprojection onto its first variable;
\item[(IV)]
  $f$ is a majority operation;
\item[(V)]
  $f$ is a minority operation.
\end{enumerate}
From now on, we will refer to an operation satisfying condition $(X)$
above, where
$X\in\{\text{I},\ \text{II}, \text{III}, \text{IV}, \text{V}\}$,
as an \emph{operation
of type $(X)$}.

The statement in the preceding paragraph can be restated as follows:
the set of all clones $\Clo(\al{A}_f)$ where
$f$ is an operation listed in~(I)--(V), 
is coinitial (with respect to inclusion) in the interval 
$(\II_A,\OO_A]$ of clones 
(containing $\OO_A$ but not containing $\II_A$).
This implies, in  particular, that
\begin{equation}\label{eq-five-pretypes}
  \begin{matrix}
  \text{every minimal clone on any set $A$ $(|A|\ge2)$ is generated by}\\
  \text{one of the operations $f$ in (I)--(V).}
  \end{matrix}
\end{equation}

If $A$ is finite, then a semiprojection on $A$ (which is not a projection,
by definition) cannot have arity $>|A|$,
and hence there are only finitely many operations of each
type~(I)--(V) on $A$. Thus, the conclusion in the preceding
paragraph implies that the analog of \eqref{eq-kuzn}
holds for minimal clones on finite sets $A$ ($|A|\ge2$):
\begin{equation}\label{eq-fin-many-min-clones}
  \begin{matrix}
    \text{for $A$ finite,
      every clone on $A$, except the clone of projections, has}\\
    \text{a minimal subclone, and the number of minimal clones
      on $A$ is finite.}
  \end{matrix}
\end{equation}  

Problem~12 in \cite[p.~120]{PK}, which we mentioned earlier, stated the task:
Describe all minimal clones on finite sets,
and determine their number.
In \cite{Ros-minclones} Rosenberg noted that analogously to the way
a description of all maximal clones on finite sets $A$ yields an efficient
criterion for checking if a set of operations on $A$ is complete
(i.e., generates the clone $\OO_A$ of all operations),
a description of all minimal clones on finite sets
$A$ would yield an efficient criterion for checking if
a relation $\rho$ on $A$ is
\emph{rigid}, that is, the polymorphism clone $\Pol(\rho)$ of $\rho$ is
the clone $\II_A$ of projections.

It is not hard to see that Statement~\eqref{eq-fin-many-min-clones} is false
if the set $A$ is infinite.
For example, if $f$ is a unary operation on $A$ which is a permutation of
infinite order, then the clone generated by $f$ has no minimal subclones.
In the second to last paragraph of Subsection~\ref{subsec-taylor-min}
we will give another example: we describe a minority operation $s$
on any infinite set $A$ such that the clone generated by $s$ has no minimal
subclones.

Let us return now to Statement~\eqref{eq-five-pretypes} saying that
every minimal clone on $A$ has a generator of one of the types~(I)--(V),
whether the base set $A$ is finite or infinite.
When we listed Examples~(i)--(viii) of minimal clones earlier, 
we described most of them by generators of types
(I)--(V); the only exceptions were the clones in Example~(ii) for $p>2$,
where the minimum arity of a non-projection is $2$.
These examples show that each type of operation listed in (I)--(V) does indeed
occur as a generator of a minimal clone.
A question that is left open by these examples is the following:
do there exist minimal clones generated by $k$-ary semiprojections on sets
of size $>k$ ($k\ge3$)? Such examples were constructed by P\'alfy in
\cite{Palfy-min-sproj-clones}. A nonconstructive proof that such
minimal clones exist was given by Lengv\'arszky~\cite{Lengvarszky}.

P\'alfy's construction in \cite{Palfy-min-sproj-clones}
yields semiprojections that satisfy strong identities
which make it fairly easy to check that the generated clones are minimal.
Rosenberg's paper~\cite{Ros-minclones} and
P\'alfy's paper~\cite{Palfy-min-sproj-clones}
are the first papers in the literature that advocated
the use of identities in classifying minimal clones.
The paper~\cite{Palfy-min-sproj-clones} also emphasized that
--- unlike the maximality of a clone ---
the minimality of a clone is an internal and abstract property.
It is internal in the sense that if 
we fix a non-projection operation $f$ on some set $A$, then
the requirement that the clone $\CC=\Clo(\al{A}_f)$ generated by $f$
be minimal 
is equivalent to the requirement that
for every non-projection operation $g\in\Clo(\al{A}_f)$, 
the algebra $\al{A}_f$ satisfies an identity of the form
\begin{equation}\label{eq-g-generates-f-back}
  \mathsf{\bold T}_2(\mathsf{\bold T}_1(f))
  (x_1,\dots,x_n)=f(x_1,\dots,x_n)
\end{equation}  
where $n$ is the arity of $f$, $x_1,\dots,x_n$ are distinct variables,
$\mathsf{\bold T}_1(f)$ is a well-formed expression representing
$g$ as a composition of copies of $f$ and projections,
and similarly,
$\mathsf{\bold T}_2(g)$ is a well-formed expression representing
$f$ as a composition of copies of $g$ and projections.
This also shows that minimality for clones is an abstract property,
that is, if $\CC$ and $\DD$ are isomorphic clones, then one of
them is minimal if and only if the other one is minimal.

In fact, more is true: if $\CC$ is a minimal clone and $\DD$ is a homomorphic
image\footnote{%
A clone homomorphism $\CC\to\DD$ is an arity preserving map
which commutes with all compositions
and maps the $i$-th $n$-ary projection in
$\CC$ to the $i$-th $n$-ary projection in $\DD$ for every $n\ge1$ and
$1\le i\le n$.}
of $\CC$ that is not the clone of projections, then $\DD$ is also
a minimal clone. In other words:
\begin{equation}\label{eq-minimality-inherited}
\begin{matrix}
  \text{If an algebra $\al{A}_f=(A;f)$ has a minimal clone, then}\\
  \text{so does every algebra $(A';f')$ in the variety
    $\var{V}(\al{A}_f)$ generated by $\al{A}_f$,}\\
  \text{whose operation $f'$ is not a projection.}
\end{matrix}
\end{equation}
If we apply~\eqref{eq-minimality-inherited}
to the $2$-element algebras
$(\mathbf{2},c_0)$, $(\mathbf{2};\neg)$, 
and $(\mathbf{2};\wedge)$ with minimal clones in Examples~(i) and (v),
then we get that every pointed set with at least $2$ elements,
every $\al{Z}_2$-set with a nontrivial action of the group $(\al{Z}_2;+)$,
and every nontrivial semilattice has a minimal clone.
Similarly, if we apply \eqref{eq-minimality-inherited} to the `ternary group'
$(\al{Z}_p;x_1-x_2+x_3)$ ($p$ prime) obtained from the additive
group $(\al{Z}_p;+)$ in Example~(ii), then we get that 
every nontrivial ternary group
$(A;x_1-x_2+x_3)$ obtained from an elementary
abelian $p$-group $(A;+)$ has a minimal clone. Since the clone of
an affine space over the field $\al{Z}_p$ is the same as the clone of
its ternary group reduct, these minimal clones are often referred to as
the clones of nontrivial affine spaces over $\al{Z}_p$.
In \cite{Plonka2}, P{\l}onka investigated 
the variety generated by an algebra $(\al{Z}_{p^2};(1-p)x_1+px_2)$ ($p$ prime)
with minimal clone in Example~(iii), found a set of defining identities
for this variety, and called this variety the
variety of \emph{$p$-cyclic groupoids}.\footnote{A \emph{groupoid} is
an algebra with one binary operation. Groupoids are also called \emph{binars}
in the literature.}
By \eqref{eq-minimality-inherited}, every groupoid
in this variety, whose operation is not a projection, has a minimal clone.

An important consequence of \eqref{eq-minimality-inherited} is
that for each type $(X)$ with
$X\in\{$I,\ II,\ III,\ IV,\ V$\}$, if $\var{W}_{(X)}$
denotes the variety in one operation symbol $f$ (of the appropriate arity)
defined
by the identities for $f$ specified in $(X)$  --- e.g.,
no identity in type (I), the identity
$f(x,x)=x$ in type (II), and
the majority identities in type (IV) ---, then the subvarieties $\var{V}$
of $\var{W}_{(X)}$ with the property that every algebra $\al{A}_f=(A;f)$
in $\var{V}$ 
with a non-projection operation $f$ has a minimal clone,
form a downward closed subset of the lattice of subvarieties of
$\var{W}_{(X)}$.
A classification of these subvarieties $\var{V}$ of $\var{W}_{(X)}$
for each type~$(X)$
would yield a classification of all
minimal clones with a fixed generating operation $f$ of type $(X)$.

This approach to classifying minimal clones
depends on a choice of the generating operation, and
there seems to be no easy way for translating between the identities if
the same clone is given in terms of another generator.
Another difficulty with this approach
is that if a clone $\Clo(\al{A}_f)$ is not minimal, i.e., if we know that
the variety $\var{V}(\al{A}_f)$ is not in the downset we are looking for, then
we might want to replace $\Clo(\al{A}_f)$ by a proper subclone
$\Clo(\al{A}_g)$ 
for some non-projection $g\in \Clo(\al{A}_f)$.
However, it may happen that $g$ is not of the same type as $f$; 
the possible type changes will be discussed at the end of the
next subsection. Furthermore, even if $g$ has the same type as $f$,
it is not clear what the relationship is between the varieties
$\var{V}(\al{A}_f)$ and $\var{V}(\al{A}_g)$.
These difficulties might be part of the reason why
the classification
of minimal clones --- even on finite sets --- has turned out to be
much more difficult than the
classification of maximal clones on finite sets.

We conclude this subsection by discussing
a technique involving identities
that has proved very useful since the early days of studying
minimal clones, and is explicitly mentioned in
\cite{Kearnes-minclone} and \cite{LevPalfy-binary-min}.
An \emph{absorption identity} in one operation symbol $f$ of arity $n$ 
is an identity of the form
\begin{equation}\label{eq-abs-id}
T_1(f)(x_1,\dots,x_m)=x_i
\end{equation}  
where $x_1,\dots,x_m$ are distinct variables and $T_1(f)$ is a well-formed
expression obtained by successive compositions from $f$ and projections.
For example, the identities defining the operations of types (II)--(V) are
the simplest absorption identities.
Absorption identities are important in the study of minimal clones
because of the following fact:
\begin{equation}\label{eq-absorp-ids}
\begin{matrix}
  \text{If $\al{A}_f=(A;f)$ has a minimal clone then $\al{A}_f$
    satisfies every}\\
  \text{absorption identity that is true in at least one algebra
  $\al{B}_f\in\var{V}(\al{A}_f)$}\\
  \text{with a non-projection basic operation.}
\end{matrix}  
\end{equation}  
Otherwise, if an absorption identity, say \eqref{eq-abs-id}, holds in
$\al{B}_f$ but not in $\al{A}_f$, then the operation
represented by the left hand side of \eqref{eq-abs-id} is a projection in
$\Clo(\al{B}_f)$ and is not a projection in 
$\Clo(\al{A}_f)$,
Therefore, the minimality of $\Clo(\al{A}_f)$ forces 
an identity of the form \eqref{eq-g-generates-f-back} in $\al{A}_f$
and hence also in $\al{B}_f$. But the latter contradicts our assumption that
the operation $f$ of $\al{B}_f$ is not a projection.

\subsection{Rosenberg's theorem on minimal clones}\label{subsec-minclone-Ros}
Let $A$ be an arbitrary set of size $\ge2$.
We know from our previous discussion that if $\CC$ is a minimal clone
on $A$, then $\CC$ is generated by an operation $f$ of one of the types
(I)--(V). For the case when $|A|=2$, say $A=\mathbf{2}$,
the converse of this statement also holds: every Boolean function $f$
of the types (I)--(V) on $\mathbf{2}$ generates a minimal clone, which appears
among the examples near the beginning of
Subsection~\ref{subsec-minclone-background}.
However, if $|A|\ge3$, then the converse statement is false.
A complete description of the minimal clones would require to identify which
operations of types~(I)--(V) generate minimal clones.
For type (I) this is fairly easy to do.
Rosenberg's theorem on minimal clones completed this step for operations of
type~(V).

\smallskip
\noindent
{\bf Rosenberg's Theorem} \cite{Ros-minclones}{\bf.}
{\it
Let $A$ be a set of size $\ge2$.
A minority operation $f$ on $A$ generates a minimal clone if and only if
$f$ is ternary addition $x_1+x_2+x_3$ for some elementary abelian $2$-group
$(A;+)$.

Consequently, every minimal clone on $A$ is generated by one of the following
types of operations $f$ on $A$:
\begin{enumerate}
\item[${}^\checkmark\!$(I)]
  $f$ is a unary operation such that either $f=f^2$ with $f(A)\not=A$, or 
  $f$ is a permutation of prime order; 
\item[(II)]
  $f$ is a binary idempotent operation that is not a projection;
\item[(III)]
  $f$ is a semiprojection onto its first variable;
\item[(IV)]
  $f$ is a majority operation;
\item[${}^\checkmark\!$(V)]
  $f$ is ternary addition $x_1+x_2+x_3$ for some elementary abelian
  $2$-group $(A;+)$. 
\end{enumerate}
The operations $f$ in {\rm ${}^\checkmark\!$(I)} and {\rm ${}^\checkmark\!$(V)}
generate minimal clones. 
}

\smallskip

The proof of the theorem in \cite{Ros-minclones}
shows that if $f$ is a minority operation
such that $f$ fails to satisfy the condition in ${}^\checkmark\!$(V),
then the clone
$\Clo(\al{A}_f)$ generated by $f$ contains a semiprojection of arity
$\ge3$.
This implies that the $2$-element subalgebras of $\al{A}_f$ satisfy absorption
identities that fail in $\al{A}_f$.  Thus,
$\Clo(\al{A}_f)$ is not a minimal clone by
\eqref{eq-absorp-ids}.

Rosenberg's Theorem is not concerned with
the question which of the operations 
of types~(II), (III), or (IV) generate minimal clones.
If one wants to follow an elimination procedure resembling the way
type~(V) was narrowed down to~${}^\checkmark\!$(V), the first step
is to answer the following question: Given an 
operation $f$ of type (II), (III), or (IV) on a set $A$,
what are the possible types, among (II)--${}^\checkmark\!$(V), for the
operations occurring in the idempotent clone
$\Clo(\al{A}_f)$?
It turns out that 
type (II) is too broad to allow any restrictions; in fact, there exist
binary idempotent operations $f$ on a $3$-element set $A$ such that
$\Clo(\al{A}_f)$ is the clone of all idempotent operations on $A$, and hence
$\Clo(\al{A}_f)$ contains all operations  of types (II)--${}^\checkmark\!$(V)
on $A$.
On the other hand,
if $f$ is of type (III) or (IV) (semiprojection or majority operation),
then we have some restrictions (see \cite{PK, Ros-minclones}).
For type~(III):
\begin{equation}\label{eq-semiproj-stable}
  \begin{matrix}\text{If $f$ is an $n$-ary semiprojection ($n\ge3$)
      onto its first variable on $A$,}\\
      \text{then every operation of one of the}\\
      \text{types (II)--${}^\checkmark\!$(V) in $\Clo(\al{A}_f)$ must
  be a semiprojection of arity $\ge n$.}
  \end{matrix}  
\end{equation}
The reason for this conclusion is that under the assumption on $f$
in~\eqref{eq-semiproj-stable}, the restriction of $f$ to every subset
$B$ of $A$ with $|B|<n$ is projection onto the first variable on $B$.
Hence, if $g\in\Clo(\al{A}_f)$, then the restriction of $g$ to every subset
$B$ of $A$ with $|B|<n$ is a projection operation on $B$.
Thus, $g$ cannot be of type (II), (IV), or
${}^\checkmark\!$(V), and $g$ cannot be a semiprojection of arity $<n$.

For type (IV):
\begin{equation}\label{eq-majority-stable}
  \begin{matrix}\text{If $f$ is a majority operation on $A$, then every
      operation of one of the}\\
    \text{types (II)--${}^\checkmark\!$(V)  in $\Clo(\al{A}_f)$ must
  be a majority operation.}
  \end{matrix}  
\end{equation}  
Indeed, if $f$ is a majority operation on $A$, then 
every operation $g\in\Clo(\al{A}_f)$ is either a projection
or a near unanimity operation
(see \cite{Csakany-conservative-minclone}).
Hence $g$ cannot be of type (II), (III), or
${}^\checkmark\!$(V).

These considerations also imply that a minimal clone cannot contain operations
of different types (I)--(V). Therefore we may unambiguously talk about a
\emph{minimal clone of type $(X)$} ($X\in\{$I, II, III, IV, V$\}$)
to mean a minimal clone generated by an operation of type $(X)$.
It will often be convenient to refer to minimal clones of types 
${}^\checkmark\!$(I), (II), (III), (IV), and ${}^\checkmark\!$(V) as
\emph{unary minimal clones}, \emph{binary minimal clones},
\emph{minimal clones of semiprojection type}, \emph{minimal majority clones},
and \emph{minimal minority clones}, respectively.

\section{Further progress on the classification of minimal clones}

As we mentioned in the preceding section, the classification of minimal clones
--- based on Rosenberg's Theorem on the five types of minimal clones ---
is far from complete for three of the five types,
even in the case when the base set is finite.
Nevertheless, over the last $40$ years or so,
significant progress has been made towards understanding these minimal clones.
In this section we will discuss some results that represent several
different directions of research on this topic, and the techniques developed
to achieve the results.\footnote{%
The reader may also be interested in the papers \cite{Quack-min-survey}
and \cite{Csakany-min-survey}
which survey results on minimal clones.}

\subsection{Minimal clones on small sets}\label{subsec-minclones-smallsets}
We saw in the preceding section that
there are seven minimal clones
on $A=\mathbf{2}$ (and hence on any $2$-element set): three 
unary, two binary, one majority, and one minority minimal clone.
This follows from Post's results in \cite{Post}, but 
can also be
derived easily from Rosenberg's Theorem (or its predecessor in
Subsection~\ref{subsec-minclone-background}).
Among these minimal clones, the two clones generated
by the constant operations $c_0$ and $c_1$,
and the two clones generated by the binary idempotent operations
$\wedge$ and $\vee$
differ only by switching the roles of the elements of the base set.
In general, we will say that two clones,
$\CC$ on $A$ and $\CC'$ on
$A'$, are \emph{similar}\footnote{%
The well-known notion of similarity for permutation groups
(see \cite[p.~32]{Robinson}) is the special case when $\CC$ and $\CC'$
are essentially unary clones whose unary parts are permutation groups.}
if there is a bijection $\pi\colon A\to A'$
such that \emph{conjugation by $\pi$}, which is a clone isomorphism
$\OO_A\to\OO_{A'}$ that sends every operation
$f(x_1,\dots,x_n)\in\OO_A$ to the operation 
${}^{\pi}f(x_1,\dots,x_n):=\pi(f(\pi^{-1}(x_1),\dots,\pi^{-1}(x_n)))\in\OO_{A'}$,
maps $\CC$ onto $\CC'$.
In this terminology, we have ${}^{\tau}c_0=c_1$ and ${}^{\tau}\wedge=\vee$
for the permutation $\tau$ switching $0$ and $1$, so on the two-element set
$A=\mathbf{2}$, the two minimal clones generated
by constant operations
and the two minimal clones generated by binary
idempotent operations are similar.
It is easy to see that no other pairs of clones among the minimal clones
on a $\mathbf{2}$ are similar, so up to similarity, there are
five minimal clones on $\mathbf{2}$: two unary, one binary, one majority,
and one minority.

The minimal clones on $3$-element sets were determined by
Cs\'ak\'any in \cite{Csakany-3el-min1} and \cite{Csakany-3el-min2}.
He found that, up to similarity, there are $24$ minimal clones on a $3$-element
set, $4$ of them are unary, $12$ of them are binary, and $3$ of them
are majority minimal clones, while the remaining 
$5$ are generated by ternary semiprojections.
The arguments were assisted by
extensive computer search, and the final result shows, for instance, that
every binary minimal clone with generating operation $f$
on a $3$-element set contains at most four binary operations: the two
projections and $f(x,y)$, $f(y,x)$ (the latter two might be equal),
which is not a priory obvious, and is not true on large
finite sets.

Now let $A$ be a $4$-element set.
The binary minimal clones on $A$ were classified by
Szczepara in his PhD thesis~\cite{Szczepara}.
He found six systems of identities (in one binary operation symbol)
such that a clone on a $4$-element set
that contains a binary non-projection is minimal if an only if it is generated
by a binary operation satisfying one of the six systems of identities.
He derived from these results that, up to similarity, 
there are exactly $120$ binary minimal clones on $A$.
The minimal majority
clones on $A$ were determined by
Waldhauser in \cite{Waldh-4el-majority-min}.
Some of these minimal clones are generated by
conservative majority operations; these minimal clones were classified
by Cs\'ak\'any~\cite{Csakany-conservative-minclone}
on any finite set, and will be discussed in the next subsection.
In \cite{Waldh-4el-majority-min},
Waldhauser proved that, up to similarity,
there are exactly $3$ minimal clones
on $A$ that are
generated by non-conservative majority operations.
For minimal clones generated by semiprojections $f$ on $A$,
the arity of $f$ can only be $3$ or $4$. 
If $f$ is $4$-ary, then $f$ is conservative, and
the classification theorem of
Je\v{z}ek and Quackenbush~\cite{JezQuack-semipr-min} --- also to be
discussed in the next subsection --- yields an efficient necessary and
sufficient condition for $f$ to generate a minimal clone on $A$.
However, if $f$ is $3$-ary, no such description is known (even in the
special case when $f$ is conservative).
Therefore,
this is the only case
missing from a complete classification of all minimal clones on $4$-element
sets. 

For finite sets $A$ of size $\ge5$, the problem of classifying all minimal
clones on $A$ is unsolved.

\subsection{Conservative minimal clones}\label{subsec-minclones-conservative}
An operation $f$ on a set $A$ is called \emph{conservative} 
if it preserves every nonempty subset $B$ of $A$, that is,
if $f\in\bigcap_{\emptyset\not=B\subseteq A}\Pol(B)$.
A clone $\CC$ is called \emph{conservative} if all operations
in $\CC$ are conservative, that is,
if $\CC\subseteq\bigcap_{\emptyset\not=B\subseteq A}\Pol(B)$.
Thus, every clone generated by conservative operations is conservative. 
It is easy to see that minimal clones of type 
${}^\checkmark$(I) in Rosenberg's Theorem are not conservative, since
conservative operations are idempotent. 
Minimal clones of type 
${}^\checkmark$(V) are not conservative
either --- unless the base set has size $2$ ---, since 
$3$-element subsets of the base set
are not closed under the generating operation
$x_1+x_2+x_3$ of the clone.
To discuss conservative minimal clones of the remaining three types, let
$f$ be a $k$-ary conservative operation on a set $A$
satisfying one of the conditions (II)--(IV) in Rosenberg's Theorem.
The assumption that $f$ is conservative is equivalent to saying that
every nonempty subset $B$ of $A$ is the base set of a subalgebra
of $\al{A}_f=(A;f)$. Here, the subalgebra of $\al{A}_f=(A;f)$
with base set $B$ is
the algebra $(B;f|_B)$ where $f|_B$ denotes the
restriction of $f$ to $B$, and hence $(B;f|_B)$ is in the variety generated by
$\al{A}_f$.
It follows that if $f$ generates
a minimal clone, then $f$ is the union of its restrictions $f|_B$
to all $k$-element subsets $B$ of $A$, and
\begin{multline}\label{eq-f-restr}
  \text{for each $f|_B$ ($B\subseteq A$, $|B|=k$), }\\
  \text{$f|_B$ generates a minimal clone on $B$ of the same
    type~(II)--(IV) as $f$,}\\
  \text{or 
    $f|_B$ is one of the two projection operation on $B$ if
    $f$ has type~(II),}\\
  \text{or 
    $f|_B$ is projection onto the first variable on $B$ if
    $f$ has type~(III).}  
\end{multline}
Thus, to classify the conservative minimal clones,
it suffices to determine under what additional conditions on the operations
$f|_B$ in \eqref{eq-f-restr}
it is the case that the union $f$ of the
$f|_B$'s generates a minimal clone on $A$.

In \cite{Csakany-conservative-minclone}, Cs\'ak\'any found such a
criterion for binary minimal clones and for minimal majority clones.
First, let $f$
be a conservative binary idempotent operation on a set $A$.
In this case, the description of all minimal clones on $2$-element sets
in Subsection~\ref{subsec-minclones-smallsets}
implies that each operation $f|_B$ in \eqref{eq-f-restr} is a
semilattice operation or one of the two projection operations on the
$2$-element set $B$.
Cs\'ak\'any proved that for a family of such operations
$f|_B$ ($B\subseteq A$, $|B|=2$), the union $f$ generates a minimal
clone on $A$
if and only if $f|_B$ is a semilattice operation for at least one $B$, and
for all projection operations $f|_B$, the algebras $(B;f|_B)$ are isomorphic
(i.e., all projection operations $f|_B$ project onto the same variable).
Now let $f$ be a majority operation on $A$.
In this case, each operation $f|_B$ in \eqref{eq-f-restr} is
a majority operation generating a minimal clone on the $3$-element set $B$.
By the description of all minimal clones on $3$-element sets
in Subsection~\ref{subsec-minclones-smallsets},
there are three minimal clones, up to similarity, which are generated by
majority operations. In one similarity class, each clone
has a unique majority operation, in another one 
each clone has three majority operations, and
in the third similarity class each clone has eight majority operations.
Cs\'ak\'any's result is that 
for a family of majority operations
$f|_B$ ($B\subseteq A$, $|B|=3$) generating minimal clones,
the union $f$ generates a minimal clone on $A$
if and only if for any two $3$-element sets $B_1,B_2\subseteq A$
such that $f|_{B_1}$ and $f|_{B_2}$ generate similar minimal clones, we have that
the algebras $(B_1,f|_{B_1})$ and $(B_2,f|_{B_2})$ are isomorphic.  

Conservative minimal clones generated by a $k$-ary semiprojection
($k\ge3$) onto the first variable
were studied by Je\v{z}ek and Quackenbush in
\cite{JezQuack-semipr-min}. Following the same strategy as in the preceding
paragraph, one has to complete the following two steps
to obtain a classification of all such minimal clones:
(a)~Describe the minimal clones on $k$-element
sets that are generated by $k$-ary semiprojections ($k\ge3$); and
(b)~Given a family of operations $f|_B$ as in \eqref{eq-f-restr}
where for each $B\subseteq A$ ($|B|=k$),
either $f|_B$ is a $k$-ary semiprojection onto the first variable that
generates a minimal clone on $B$, or
$f|_B$ is the $k$-ary projection onto the first variable on $B$,
find a necessary and sufficient condition for the union $f$ of these $f|_B$'s
to generate a minimal clone on $A$.
To state the results of \cite{JezQuack-semipr-min} on Steps~(a) and (b), we
call a binary relation $\rho$ on a set $B$ \emph{bitransitive} if
$\rho$ is reflexive and transitive, but not the equality relation,
and the automorphism group
of the relational structure $(B;\rho)$ acts transitively on the set
of all pairs $(b_1,b_2)\in\rho$ with $b_1\not=b_2$.
Bitransitive relations on $B$ fall into two types: 
equivalence relations on $B$ with all non-singleton
blocks of the same size, and partial orders of length $1$ whose
automorphism group acts transitively on the set of pairs $(b_1,b_2)\in\rho$
with $b_1\not=b_2$.

Je\v{z}ek and Quackenbush completed Step~(a) 
by proving the following:
If $B$ is a $k$-element set and $\CC$ is a clone on $B$ generated by a $k$-ary
semiprojection, then $\CC$ is minimal if and only if $\CC$ is generated by
a semiprojection of the form
\begin{equation*}
  s_\rho(b_1,\dots,b_k)=
  \begin{cases} b_k & \text{if $b_1,\dots,b_k$ are distinct
      and $(b_1,b_k)\in\rho$,}\\
    b_1 & \text{otherwise}
  \end{cases}
\end{equation*}
for some bitransitive relation $\rho$ on $B$.
Moreover, for any two distinct bitransitive relations $\rho_1$ and $\rho_2$
on $B$,
the operations $s_{\rho_1}$ and $s_{\rho_2}$ generate distinct minimal clones.
For Step~(b), Je\v{z}ek and Quackenbush proved that the following condition
is necessary 
for a conservative clone $\CC$ generated by
a $k$-ary semiprojection on a set $A$ to be minimal:
$\CC$ is generated by a semiprojection $f$ such that
\begin{enumerate}
\item[$(*)$]
  all restrictions $f|_B$ ($B\subseteq A$, $|B|=k$)
  of $f$ that are not projections have the form $f|_B=s_{\rho_B}$
  for isomorphic bitransitive relational structures $(B;\rho_B)$.
\end{enumerate}  
This necessary condition for $\CC$ to be minimal is not sufficient if $|A|>k$.
If $f$ generates $\CC$ and satisfies $(*)$, then the algebra
$\al{A}_f=(A;f)$ and
its subalgebras $(B;f|_B)$, where $f|_B$ is not a projection, will satisfy
the same $k$-variable identities, but the clone $\CC=\Clo(\al{A}_f)$
may contain a semiprojection of arity $>k$, which implies
by \eqref{eq-absorp-ids} that $\CC$ is not minimal.
Therefore,
by adding the condition ``$(A;f)$ satisfies the same absorption identities
as its (isomorphic) subalgebras $(B;f|_B)$ (with $B$ as in $(*)$)
that are not projection algebras''
to the necessary condition above, one can get a characterization of
all conservative minimal clones of semiprojection type.
However, as it was pointed out
in \cite{JezQuack-semipr-min}, there seems to be no efficient way for checking
this last condition.

\subsection{Minimal clones with few non-projection operations
  of minimum arity}\label{subsec-minclone-few-ops}
The binary minimal clones that we have encountered so far
don't have `too many'
binary non-projection operations: in Examples~(ii)--(v) in
Subsection~\ref{subsec-minclone-background}
the numbers are $p-2$, $2p-2$, $2$, and $1$
respectively, while
for the binary minimal clones on $2$- or $3$-element sets
and for the conservative
binary minimal clones the numbers are $1$ or $2$
(see Subsections~\ref{subsec-minclones-smallsets} and
\ref{subsec-minclones-conservative}).
In the mid 1990's, 
L\'evai and P\'alfy~\cite{LevPalfy-binary-min}
started a systematic study of binary minimal clones
$\CC$ in which the number of binary non-projections,
$|\CC^{(2)}|-2 (>1)$, is small.
Their results classify all binary minimal clones $\CC$ with
$|\CC^{(2)}|=3, 4, 6$.
When combined with two theorems of Dudek
for the cases $|\CC^{(2)}|=5, 7$
--- one in~\cite{Dudek} and the other unpublished at the time,
but published much later in~\cite{Dudek-Galuszka} ---
these results yield a classification of all binary minimal clones
$\CC$ with $|\CC^{(2)}|\le 7$.

Since an assumption that restricts the size of the binary part
$\CC^{(2)}$ of a binary minimal clone $\CC$ 
depends only on the isomorphism type of $\CC$,
these classification results determine the minimal clones $\CC$
considered up to isomorphism.
The classification is accomplished by finding,
under each assumption $|\CC^{(2)}|=r$ ($r=3,\dots,7$),
a family
$\mathcal{F}_r$ of groupoid varieties
--- i.e., varieties in one
binary operation symbol $\circ$ ---
such that
\begin{equation}\label{eq-LP}
  \begin{matrix}
    \text{a clone $\CC$ on a set $A$ is a binary minimal
      clone satisfying $|\CC^{(2)}|=r$}\\
    \text{if and only if $\CC$ contains a binary
      non-projection $\circ$ such that}\\
    \text{$\al{A}_\circ=(A;\circ)$ is in $\var{V}$
      for some $\var{V}\in\mathcal{F}_r$,
      $\CC=\Clo(\al{A}_\circ)$, and $|\CC^{(2)}|\not<r$.}
  \end{matrix}  
\end{equation}  
For $r=3$ there are two varieties in $\mathcal{F}_3$:
\begin{itemize}
\item
  the variety $\var{A}_3$ of affine spaces over $\al{Z}_3$;
\item
  the variety $2$-$\mathcal{SL}$
  defined by the $2$-variable identities of the variety
  $\mathcal{SL}$ of semilattices.
\end{itemize}
For $r=4$, $\mathcal{F}_4$ consists of five varieties,
which we also list here\footnote{In some cases, our description
of these varieties is different from, but equivalent to, 
the description in \cite{LevPalfy-binary-min}.
The terminology `meld' for the members of the last variety 
is due to Brady~\cite{Brady-binary-min-nontaylor}.}
because they will play a role later on:
\begin{itemize}
\item
  the variety $2$-$\mathcal{CG}$ of $2$-cyclic groupoids;
\item
  the variety $\mathcal{RB}$ of rectangular bands;
\item
  the variety $2$-$\mathcal{LNB}$
  defined by the $2$-variable identities of the variety
  $\mathcal{LNB}$ of left normal bands,
  where $\mathcal{LNB}=\mathcal{SL}\vee\mathcal{LZ}$ (the
  join of the variety of semilattices and the variety of
  left zero semigroups);
\item
  the variety $\mathcal{D}$ defined by the identities
  \begin{gather}
    x\circ(y\circ x)=(x\circ y)\circ x=
    (x\circ y)\circ y=x\circ y,\quad\text{and}\label{eq-LP-d1}\\
    x\circ((\dots((x\circ y_1)\circ y_2)\dots)\circ y_m)=x
    \quad \text{for all }m=0,1,2\dots;
    \label{eq-LP-d2}
  \end{gather}  
\item
  the \emph{variety $\var{M}$ of melds}, which is defined by the identities
  \begin{equation*}
    x\circ x=x,\quad  x\circ((y\circ x)\circ z)=x.
  \end{equation*}  
\end{itemize}  
For $r=5$ and
$r=7$, $\mathcal{F}_r$ is a singleton containing the variety of 
affine spaces over $\al{Z}_r$. Finally, for $r=6$,
$\mathcal{F}_6$ consists of five varieties, including the variety
of $3$-cyclic groupoids.

The method L\'evai and P\'alfy used in \cite{LevPalfy-binary-min} to find
the varieties in $\mathcal{F}_r$ for $r=3,4,6$ is based on the idea that
if $\CC=\Clo(\al{A}_\circ)$ is a binary minimal clone
with $|\CC^{(2)}|=r$, then
the $r$-element groupoid $\al{C}_\circ:=(\CC^{(2)};\circ)$,
where $\circ$ is the binary operation
\begin{equation*}
  \CC^{(2)}\times\CC^{(2)}\to\CC^{(2)},\quad
  \bigl(f(x_1,x_2),g(x_1,x_2)\bigr)\mapsto f(x_1,x_2)\circ g(x_1,x_2)
\end{equation*}
of applying the operation $\circ$
of $\al{A}_\circ$ pointwise to the functions in $\CC^{(2)}\,(\subseteq A^{A^2})$,
has two important properties. First,
\begin{enumerate}
\item[$(\alpha)$]
  $\al{C}_\circ$ is the $2$-generated free groupoid in the variety
  $\var{V}(\al{A}_\circ)$, with the two projections
  $\pr_1^{(2)},\pr_2^{(2)}\in\CC^{(2)}$
  as free generators.
\end{enumerate}
  This implies that
  the binary parts of the clones $\CC$ and $\Clo(\al{C}_\circ)$
  are isomorphic. Therefore, since $\CC$ is a minimal clone, we get the second
  property that
\begin{enumerate}
\item[$(\beta)$]
  the clone of $\al{C}_\circ$ is \emph{$2$-minimal}, that is,
  every binary non-projection in $\Clo(\al{C}_\circ)$ generates
  the operation $\circ$ of $\al{C}_\circ$.
\end{enumerate}  
Properties $(\alpha)$ and $(\beta)$
narrow down the isomorphism type of the groupoid $\al{C}_\circ$
to $2$, $7$, and $6$ possibilities for
$r=3$, $4$, and $6$, respectively, provided groupoids with the same clone
are counted only with one choice of the generating operation
$\circ$. (In the case $r=6$, the search was done by computer.)
For some of these isomorphism types, $\al{C}_\circ$
does not have a minimal clone.
These isomorphism types can be eliminated, because
we have that $\al{C}_\circ\in\var{V}(\al{A}_\circ)$
and the operation of $\al{C}_\circ$
is a non-projection, hence \eqref{eq-minimality-inherited} forces
the clone $\Clo(\al{C}_\circ)$ to be minimal.
Finally, for each $r=3,4,6$ and
each isomorphism type, say $i=1,\dots,k_r$, for an $r$-element
$\al{C}_\circ$ with a minimal clone,  
the authors of~\cite{LevPalfy-binary-min}
found a set of identities, $\Sigma_{r,i}$, consisting of binary identities and 
absorption identities true in $\al{C}_\circ$ such that whenever
a groupoid with a non-projection operation
satisfies the identities in $\Sigma_{r,i}$,
then the groupoid has a minimal clone.
By the construction of $\al{C}_\circ$ and by 
Statement~\eqref{eq-absorp-ids}, this implies that
\eqref{eq-LP} holds if for each $r=3,\,4,\,6$, $\mathcal{F}_r$
is chosen to be the set of varieties $\var{V}_{r,i}$ 
defined by the identities in $\Sigma_{r,i}$ ($i=1,\dots,k_r$).

In \cite{LevPalfy-binary-min},
L\'evai and P\'alfy also constructed
binary minimal clones $\BB_k$ and $\CC_k$, for each $k\ge2$, such that
$|\BB_k^{\,(2)}|=2k+2$ and $|\CC_k^{\,(2)}|=3k+2$. These examples show that
there is no upper bound on the number of binary non-projections in binary
minimal clones. The second sequence also shows that for infinitely many primes
$p\,(\ge 11)$ of the form $3k+2$, the clones of nontrivial affine spaces over
$\al{Z}_p$ are not the only minimal clones 
$\CC$ with $|\CC^{(2)}|=p$.

Now we will discuss minimal majority clones with few ternary operations.
Recall from \eqref{eq-majority-stable}
that in a minimal majority clone, every ternary non-projection is a
majority operation.
The known examples of minimal majority clones $\CC$ exhibit very few
possibilities for the number $|\CC^{(3)}|-3$ of majority operations in $\CC$.
For Examples~(vi) and~(vii) in
Subsection~\ref{subsec-minclone-background},
the numbers are $1$ and $3$, respectively. 
For the majority minimal clones on $3$-element sets
and for the non-conservative majority minimal clones on $4$-element sets
the numbers are $1$, $3$, or $8$
(see Subsection~\ref{subsec-minclones-smallsets}).
For the conservative
minimal majority clones the numbers are $1$, $3$, $8$, or $24$
(see Subsection~\ref{subsec-minclones-conservative}).
In the paper \cite{Waldh-majority-min}, Waldhauser classified all
minimal majority clones with $\le7$ ternary operations (i.e., with $\le4$
majority operations).
The main results are that there exist no minimal majority clones
with $2$ or $4$ majority operations, and for each one of $m=1$ and $3$,
there exists a unique variety $\var{V}_m$ (in one ternary operation
symbol $f$) such that
\begin{equation}\label{eq-Waldh}
  \begin{matrix}
    \text{a clone $\CC$ on a set $A$ is a minimal majority clone
      satisfying $|\CC^{(3)}|=3+m$}\\
    \text{(i.e., $\CC$ contains exactly $m$ majority operations)}\\
    \text{if and only if $\CC$ contains a majority operation
      $f$ such that}\\
    \text{$\al{A}_f\in\var{V}_m$,\ \ 
      $\CC=\Clo(\al{A}_f)$,\ \ and\ \ $|\CC^{(3)}|\not<m$.}
  \end{matrix}  
\end{equation}  
Here, $\var{V}_1$ is the variety defined by the $3$-variable identities that
are true for the median operation of any lattice with $2$ or more elements,
while $\var{V}_3$ is the variety defined by the $3$-variable identities
that are true for the dual discriminator operation on any set of size $\ge3$
(see Examples~(vi) and~(vii) in Subsection~\ref{subsec-minclone-background}).
Clearly, $\var{V}_1\subset\var{V}_3$.
Waldhauser also pointed out in \cite{Waldh-majority-min} that $\var{V}_1$
has infinitely many subvarieties, and hence there are infinitely many
non-isomorphic minimal majority clones that contain a single majority operation.

In the classification of minimal majority clones $\CC$ by identities,
it is the ternary part, $\CC^{(3)}$, that plays the same role
as the binary part for binary minimal clones.
However, there is an essential difference: for any clone $\CC=\Clo(\al{A}_f)$
generated by a majority operation $f$, 
Statement \eqref{eq-majority-stable} implies that
$\CC$ is minimal if (and only if) $\CC$ is 
\emph{$3$-minimal}, that is,
every ternary non-projection $g$ in $\CC$ generates $f$.
Thus, the ternary part $\CC^{(3)}$ determines whether or not $\CC$ is minimal,
or equivalently, the ternary identities true in $\al{A}_f$
determine whether or not the clone
$\Clo(\al{A}_f)$ is minimal.

The main tool Waldhauser used in his classification of minimal majority clones
$\CC$ with $|\CC^{(3)}|\le7$ is the symmetries of the majority operations in
$\CC$. The symmetric group $S_3$ acts on
$\CC^{(3)}$ by permuting the variables of the operations.
Call an operation $f\in\CC^{(3)}$
\emph{cyclically symmetric} if it is invariant under permuting its variables
by a $3$-cycle in $S_3$.
Waldhauser proved for every minimal majority clone $\CC$
with finitely many majority operations that
if every majority operation in
$\CC$ is cyclically symmetric, then $\CC$ has a unique majority operation.
Consequently, if $\CC$ is a minimal majority clone
with finitely many, but more than one majority operations, then
some majority operation in $\CC$ has an orbit of size $\ge3$ under the action
of $S_3$. This shows that the number $m$ of majority operations in
a minimal majority clone $\CC$ cannot be $2$, and if $m=3$ then the
three majority operations in $\CC$ are obtained from one another
via permuting the variables by $3$-cycles. Finally, were $m=4$, there would be
one cyclically symmetric majority operation in $\CC$ and three other majority
operations that are obtained from one another via permuting variables
by $3$-cycles. This
is the starting point of eliminating the case $m=4$.

For a long time, the only numbers known to occur as the number of majority
operations in a minimal majority clone were the numbers
$m=1,\,3,\,8$, and $24$ mentioned above.
However, in their paper \cite{BehrWald-majority-min} published in
2011, Behrisch and Waldhauser found minimal majority clones
$\MM_{26}$ and $\MM_{78}$ with $m=26$ and $78$ majority operations, respectively.
The clone $\MM_{26}$ is a minimal majority clone on a $5$-element set, generated
by a cyclically symmetric majority operation. It was found by a
computer search that took several weeks on several computers to
find all minimal majority clones generated by cyclically symmetric
majority operations on a $5$-element set,
and determine the number of majority operations in them.
The search revealed that, up to similarity, $\MM_{26}$ is the only such clone
where the number of majority operations is different from $1$ and $8$.
The clone $\MM_{78}$ was obtained from $\MM_{26}$ by a general construction
in the paper \cite{BehrWald-majority-min}, which produces from any
minimal majority clone with $m$ majority operations, which is generated
by a cyclically symmetric majority operation $f$ on a set $A$, a new
minimal majority clone with $3m$ majority operations, which is generated
by the majority operation $f^*$ on the set $A^*=A\cup\{*\}$ that extends
$f$ so that $f^*(a_1,a_2,a_3)=a_1$ whenever $\{a_1,a_2,a_3\}$ is a $3$-element
subset of $A^*$ containing the new element~$*$.

\subsection{Minimal clones satisfying some algebraic
  restrictions}\label{subsec-minclone-alg-restr}
Abelian algebras were introduced earlier in Subsection~\ref{subsec-abelian},
and it was mentioned there that
typical examples are the algebras whose clone is a
subclone of the clone of the constant expansion of a module,
and their subalgebras.
In \cite{Kearnes-minclone}, Kearnes classified all minimal clones
that have nontrivial abelian representations.
Given a minimal clone $\CC$ on some set $S$, we call an algebra
$\al{A}_f=(A;f)$ a \emph{representation of $\CC$}
if there exists a surjective clone homomorphism
$\Phi\colon\CC\to\Clo(\al{A}_f)$, or equivalently,
if
for an appropriately chosen non-projection\footnote{%
$\tilde{f}\in\CC$ is one of the inverse images of $f$ under $\Phi$.}
$\tilde{f}\in\CC$ 
the algebra $\al{S}_{\tilde{f}}=(S;\tilde{f})$, whose clone is $\CC$,
generates a variety $\var{V}(\al{S}_{\tilde{f}})$ containing $\al{A}_f$.
A representation $\al{A}_f$ of $\CC$ is called \emph{faithful} if
there exists a clone isomorphism $\Phi\colon\CC\to\Clo(\al{A}_f)$, or
in the equivalent reformulation, $\al{A}_f$ is another generating algebra
for $\var{V}(\al{S}_{\tilde{f}})$.
We say that a representation $\al{A}_f$ of $\CC$ is \emph{abelian}
if $\al{A}_f$ is an abelian algebra, and \emph{nontrivial} if $f$
is a non-projection.
Recall that by Statement~\eqref{eq-minimality-inherited},
a nontrivial representation $\al{A}_f$ of a minimal clone $\CC$ has a minimal
clone. 
The descriptions of the minimal clones in Examples (i)--(iv)
in Subsection~\ref{subsec-minclone-background} show that these
are clones of abelian algebras, and hence these clones have
nontrivial abelian representations (where $\Phi$ is the identity map).

The main result of \cite{Kearnes-minclone} is that, up to isomorphism,
these examples include almost all of the minimal clones with abelian
representations. Using the terminology introduced right after
Statement~\eqref{eq-minimality-inherited}, Kearnes' theorem
can be stated as follows:

\smallskip
\noindent
{\it $\CC$ is a minimal clone with
a nontrivial abelian representation if and only if one of the following
conditions holds for $\CC$:
\begin{enumerate}
\item[(o)]
  $\CC$ is unary,
\item[(a)]
  $\CC$ is the clone of an affine space over
  the field $\al{Z}_p$ for some prime $p$,
\item[(b)]
  $\CC$ is the clone of a rectangular band which is not a left zero
  or right zero semigroup, or
\item[(c)]
  $\CC$ is the clone of a
  $p$-cyclic groupoid for some prime $p$ which is not a left zero semigroup.
\end{enumerate}}

\smallskip
The proof starts by observing that among the five types of
the minimal clones (see Rosenberg's Theorem in
Subsection~\ref{subsec-minclone-Ros}) those of
types~${}^\checkmark$(I) and ${}^\checkmark$(V) are clones of abelian algebras,
and those of types~(III) and (IV) cannot have nontrivial abelian
representations. Therefore, most of the work in the proof
goes into analyzing the minimal clones of
type~(II), that is, the binary minimal clones. In the first phase of the
proof Kearnes shows that every binary minimal clone that has a nontrivial
abelian representation, is represented by one of the abelian
algebras $\al{A}_f=(\al{Z}_6;3x_1+4x_2)$,
$\al{A}_f=(\al{Z}_p;x_1-x_2+x_3)$ ($p>2$ prime),
or $\al{A}_f=(\al{Z}_{p^2};(1-p)x_1+px_2)$ ($p$ prime),
whose clones are exactly the minimal clones
in Examples (iv), (ii), and
(iii) in Subsection~\ref{subsec-minclone-background}.
In the second phase he shows that each variety generated by one of
these algebras $\al{A}_f$
can be axiomatized by absorption identities, and hence 
Statement~\eqref{eq-absorp-ids} implies that
the minimal clone $\CC=\Clo(\al{S}_{\tilde{f}})$
represented by $\al{A}_f$ must arise from an algebra $\al{S}_{\tilde{f}}$
that satisfies
all identities of $\al{A}_f$, implying that $\al{A}_f$ is a faithful
representation of $\CC=\Clo(\al{S}_{\tilde{f}})$.

The classification of minimal clones
with nontrivial abelian representations, in combination with some results from
commutator theory, 
allowed the author of \cite{Kearnes-minclone} to prove the following
generalization of the main statement in Rosenberg's Theorem on
minority minimal clones:

\smallskip
\noindent
{\it
If a clone $\CC$ is generated by a Mal'tsev operation (i.e., by a ternary
operation $p$ satisfying the identities in \eqref{eq-maltsev}),
then $\CC$ has a nontrivial abelian representation, and hence 
$\CC$ is the clone of an affine
space over a field $\al{Z}_p$ for some prime $p$.}

\smallskip
\noindent
In \cite{Waldh-minclones-weakly-abelian}, Waldhauser proved that
Kearnes' classification also holds for the minimal clones
that have weakly abelian representations.

It is not hard to see from the classification of minimal clones with abelian
representations that all these clones are commutative.
A clone $\CC$ on a set $A$ is called \emph{commutative} if any two operations
in $\CC$ \emph{commute} with each other, that is, for any $f,g\in\CC$ (say 
$f$ is $m$-ary and $g$ is $n$-ary), $g$ is an algebra homomorphism
$(A;f)^m\to(A;f)$, or equivalently, $f$ preserves the graph of $g$
(an $(m+1)$-ary relation). Yet another equivalent condition is that $f$ and $g$
satisfy the following identity in $mn$ variables $x_{ij}$:
\begin{multline}\label{eq-commutativity}
  g(f(x_{11},\dots,x_{m1}),\dots,f(x_{1n},\dots,x_{mn}))\\
  =f(g(x_{11},\dots,x_{1n}),\dots,g(x_{m1},\dots,x_{mn})).
\end{multline}  
It follows from these equivalent descriptions of commutativity
that $f$ and $g$ have symmetric roles in the definition.
It follows also that
the commutativity of clones is
invariant under clone isomorphism, and 
if $\CC$ is a minimal clone
generated by $f$, then $\CC$ is a commutative clone if and only if $f$
commutes with itself. These facts are sufficient to derive from
the classification of minimal clones with abelian
representations that all such minimal clones are commutative.

This observation motivated the classification of all commutative minimal
clones in \cite{KSz-commutative-min} by Kearnes and Szendrei.
By our discussion in the preceding paragraph,
the minimal clones of types
${}^\checkmark$(I) and ${}^\checkmark$(V) are commutative.
On the other hand, it is not hard to
show that no majority operation on a set of size $>1$ commutes with itself,
so minimal clones of types (IV) are not commutative.
Therefore the main task in \cite{KSz-commutative-min} was to classify 
the commutative minimal clones of types~(II) and (III).
We will start with the binary case.
If $\circ$ is a binary operation, then the identity expressing that $\circ$
commutes with itself is the $4$-variable identity
$(x\circ y)\circ(u\circ v)=(x\circ u)\circ(y\circ v)$,
which is known in the literature as the \emph{entropic law}
(or \emph{medial law}).
Therefore, to classify the commutative binary minimal clones, one has to
classify the idempotent entropic groupoids $\al{A}_\circ=(A;\circ)$
with minimal clones.
The main result of \cite{KSz-commutative-min}
for this case is that $\CC$ is a commutative
binary minimal clone if and only
if $\CC$ is one of the clones in items (a)--(c) above, or
\begin{enumerate}
\item[(d)]
  $\CC$ is the clone of a \emph{right semilattice}
  --- i.e., an idempotent entropic groupoid $\al{A}_\circ$
  satisfying the identities $x\circ(x\circ y)=x$ and
  $(x\circ y)\circ y=x\circ y$ ---
  which is not a left zero semigroup, or
\item[(e)]
  $\CC$ is the clone of a left normal band
  which is not a left zero semigroup.
\end{enumerate}
In both cases (d) and (e) the binary part,
$\CC^{(2)}$, of the clone $\CC$ has size $4$ or $3$ (the latter only
in case (e) when $\CC$ is the clone of a semilattice),
therefore for a non-projection $\circ$ in $\CC$, the algebra $\al{A}_\circ$
belongs to one of the varieties in $\mathcal{F}_4$ or
$\mathcal{F}_3$ in the classification theorem of L\'evai and P\'alfy in 
\cite{LevPalfy-binary-min} (see Subsection~\ref{subsec-minclone-few-ops}).
In case~(e) the variety is clearly
$2$-$\mathcal{LNB}\in\mathcal{F}_4$ or
$2$-$\mathcal{SL}\in\mathcal{F}_3$, respectively, while in
case~(d) it is not hard to check that the variety is
$\mathcal{D}\in\mathcal{F}_4$.

In contrast to the binary case, the classification in 
\cite{KSz-commutative-min} of the commutative
minimal clones generated by semiprojections produced
a new class of minimal clones.
Indeed, for minimal clones of type~(III),
we have seen so far only classification results for
minimal clones generated by conservative semiprojections.
However, it is not hard to show
that such a clone cannot be commutative.

To present the result of \cite{KSz-commutative-min} for the
commutative minimal clones generated by semiprojections,
let $\CC$ be a clone on a set $A$ generated by a $(k+1)$-ary semiprojection
$s=s(x,\wec{y})$
onto the variable $x$ (where $k\ge2$ and
$\wec{y}=(y_1,\dots,y_k)$ is the tuple of
the remaining variables), and assume that $\CC$ is commutative,
that is, $s$ commutes with itself.
The assumption that $s$ is a semiprojection and generates $\CC$ implies that
every operation $g=g(x_1,\dots,x_n)\in\CC$ restricts to all $k$-element
subsets of $A$ as projection onto the same variable $x_\ell$, which we will call
the \emph{distinguished variable} of $g$. Hence, in particular,
every operation in $\CC^{(k+1)}$ is a semiprojection or a projection.
The assumption that $s$ commutes with itself, that is, 
the assumption that the
identity~\eqref{eq-commutativity} holds with $f=g=s$ and $m=n=k+1$,
implies that the $(2k+1)$-variable
identity $s(s(x,\wec{y}),\wec{z})=s(s(x,\wec{z}),\wec{y})$ also holds for
$s$, and
every operation $t\in\CC^{(k+1)}$ with distinguished variable $x$
and non-distinguished variables $\wec{y}=(y_1,\dots,y_k)$ has the form
$s(\dots s(s(x,\wec{y}_{\sigma_1}),\wec{y}_{\sigma_2}),\dots,\wec{y}_{\sigma_m})$
for some $m\ge1$ and some permutations $\sigma_1,\dots,\sigma_m\in S_k$,
where $\wec{y}_\sigma$ denotes the tuple $(y_{\sigma(1)},\dots,y_{\sigma(k)})$
for every permutation $\sigma\in S_k$.
Moreover, every subclone of $\CC$, except the clone $\II_A$
of projections,
contains a $(k+1)$-ary semiprojection.
Thus, $\CC$ is completely determined by its subset $M_{\CC}$ consisting
of all $(k+1)$-ary operations (semiprojection or projection operations)
with a fixed distinguished variable,
say $x$. The properties of $\CC$ just discussed imply that 
the binary operation $\oplus$ on $M_{\CC}$ defined by
$(t_1\oplus t_2)(x,\wec{y}):=t_2(t_1(x,\wec{y}),\wec{y})$
yields a commutative monoid whose neutral element  
$o\in M_{\CC}$ is the $(k+1)$-ary projection onto the variable $x$.
Furthermore, the unary operations
$t(x,\wec{y})\mapsto (\sigma t)(x,\wec{y}):=t(x,\wec{y}_\sigma)$
($\sigma\in S_k$)
on $M_{\CC}$ define an action of $S_k$ on the monoid $(M_\CC;\oplus,o)$
by automorphisms, yielding an $S_k$-semimodule
$\al{M}_\CC=(M_\CC;\oplus,o,S_k)$.

It is proved in \cite{KSz-commutative-min} that
the assignment $\CC\mapsto\al{M}_\CC$ is the object map of a categorical
equivalence between the category of commutative clones generated by
a $(k+1)$-ary semiprojection or projection,
together with all clone homomorphisms between them, and
the category of $1$-generated $S_k$-semimodules with all homomorphisms
between them. This implies, in particular, that $\CC$ is a commutative
minimal clone if and only if its associated $S_k$-semimodule $\al{M}_\CC$
is \emph{minimal} in the sense that it is nontrivial (i.e., $M_\CC\not=\{o\}$)
and every element $f\in M_\CC\setminus\{o\}$ generates $\al{M}_\CC$.
The minimal $S_k$-semimodules are classified in
\cite{KSz-commutative-min} as follows:
A minimal $S_k$-semimodule is either a $2$-element semilattice
$(\{s,o\};\oplus,o)$ with neutral element $o$ and with the trivial action
of $S_k$, or an elementary abelian $p$-group $(M;\oplus,o)$
for some prime $p$ with an action
of $S_k$ that makes it an irreducible $S_k$-module over the field $\al{Z}_p$.

This yields a classification of all commutative minimal clones generated by
a semiprojection as follows.
Since there is a bijection between the isomorphism classes of
commutative minimal clones generated by $(k+1)$-ary semiprojections
and the isomorphism
classes of 
minimal $S_k$-semimodules, it suffices to describe how to construct
from each minimal semimodule $\al{M}=(M;\oplus,o,S_k)$
a commutative minimal clone $\CC$ generated by a
$(k+1)$-ary semiprojection such that
$\al{M}\cong\al{M}_\CC$.
Given a minimal semimodule $\al{M}=(M;\oplus,o,S_k)$, 
choose any $s_0\in M\setminus\{o\}$ and fix a presentation
of $\al{M}$ with generator $s_0$. The relations of this presentation
yield a set of identities in one $(k+1)$-ary operation symbol $s$ which,
together with the semiprojection identities for $s$
and the identity expressing that $s$ commutes with itself,
define a variety $\var{V}_{\al{M}}$ such that $\al{M}\cong\al{M}_\CC$ holds
for the clone $\CC=\Clo(\al{A}_s)$ of each algebra $\al{A}_s=(A;s)$ 
generating $\var{V}_{\al{M}}$. In fact, $\al{A}_s$ will generate
$\var{V}_{\al{M}}$ as long as the operation $s$ is not a projection,
as we will now show.
Since $\al{M}\cong\al{M}_\CC$ is a minimal $S_k$-semimodule, $\al{M}_\CC$
is a simple algebra, that is,
the only proper homomorphic image of $\al{M}_\CC$ is the trivial (one-element)
$S_k$-semimodule. This property carries over to $\CC$, since
$\CC$ and $\al{M}_\CC$ correspond to each other under
the categorical equivalence discussed in the preceding paragraph. Hence,
the only proper homomorphic image of
$\CC=\Clo(\al{A}_s)$, for an algebra $\al{A}_s$ generating $\var{V}_{\al{M}}$,
is the clone of projections.
But this fact is equivalent to the fact that 
the only nontrivial proper subvariety of
$\var{V}_{\al{M}}$ is the variety of algebras $\al{A}_s$ where $s$
is a projection.

Altogether this shows that commutative minimal clones generated
by semiprojections are classified by irreducible representations
of finite symmetric groups over fields of prime order.

\subsection{Brady's coarse classification of the binary minimal
  clones on finite sets}\label{subsec-brady}  
We discussed in Subsections~\ref{subsec-minclone-Ros} and
\ref{subsec-minclone-few-ops}
that one of the difficulties with classifying binary minimal clones is
that a clone $\CC$ generated by a binary idempotent non-projection may have
minimal subclones that don't have non-projection binary operations, and this may
happen even if $\CC$ is $2$-minimal. This is in sharp contrast to the case
when $\MM$ is generated by a majority operation, because
this family of clones is `stable' under taking minimal subclones in the sense
that every minimal subclone of a clone generated by a majority operation is
generated by a majority operation, and hence $3$-minimality implies minimality
for clones generated by majority operations.

In his recent papers~\cite{Brady-binary-min-taylor} and
\cite{Brady-binary-min-nontaylor}, Brady found a coarse
classification for binary minimal clones on finite sets which eliminates the
issue of `instability' mentioned in the preceding paragraph.
More precisely, Brady's main result isolates seven disjoint families of
clones with the following properties:
(1)~every clone in each family is generated by a binary idempotent
non-projection on a finite set,
(2)~every binary minimal clone on a finite set belongs to one of
the seven families,
(3)~four of the seven families consist of binary minimal clones only,
and
(4) for each one of the remaining three
families, if $\CC$ is in a given family, then every minimal subclone of $\CC$
is in the same family.
The classification is as follows: 

\smallskip
\noindent
{\it Every binary minimal clone $\CC$ on a finite set satisfies
exactly one of the following seven conditions:
\begin{enumerate}
\item[(a)]
  $\CC$ is the clone of an \ul{a}ffine space over the field $\al{Z}_p$
  for some odd prime $p$;
\item[(b)]
  $\CC$ is the clone of a rectangular \ul{b}and which is not a left zero or
  right zero semigroup;
\item[(c)]
  $\CC$ is the clone of a $p$-\ul{c}yclic groupoid for some prime $p$
  which is not a left zero semigroup;
\item[(m)]
  $\CC$ is the clone of a \emph{\ul{m}eld}
  (see Subsection~\ref{subsec-minclone-few-ops}) which is
  not a left zero semigroup;
\item[(B)] 
  $\CC$ is the clone of a \emph{\ul{B}rady groupoid}\footnote{%
In \cite{Brady-binary-min-nontaylor},
Brady called these groupoids `partial semilattices'.
Since we discuss partial operations and partial clones
in Subsections~\ref{subsec-max-partial-clones}
and~\ref{subsec-min-partial-clones} of this paper,
it seemed safer to avoid this name.}  ---
  i.e., a groupoid $\al{A}_\circ$
  satisfying the identities $x\circ x=x$ and 
  $x\circ (x\circ y) = (x\circ y)\circ x = x\circ y$ --- which is
  not a left zero semigroup, but has a $2$-element section (i.e., quotient
  of a subalgebra) that is a
  left zero semigroup;
\item[(S)]
  $\CC$ is the clone of a \emph{\ul{s}piral} ---
  that is, a groupoid $\al{A}_\circ$
  satisfying the identities $x\circ x=x$, $x\circ y=y\circ x$, and having
  the additional property that
  for any two distinct elements $a,b\in A$, the
  subgroupoid of $\al{A}_\circ$ generated by $\{a,b\}$
  either contains only these two
  elements, or has a surjective homomorphism
  onto the free semilattice $\al{F}_{\mathcal{SL}}(x,y)$ on two generators;
\item[(D)]
  $\CC$ is the clone of a \emph{\ul{d}ispersive groupoid} --- that is,
  a groupoid $\al{A}_\circ$
  satisfying all identities in \eqref{eq-LP-d2} and having the additional
  property that
  for any two distinct elements $a,b\in A$, either the
  subgroupoid of $\al{A}_\circ$ generated by $\{a,b\}$ contains only these two
  elements, or
  the subgroupoid of $(\al{A}_\circ)^2$ generated by the pairs
  $(a,b)$ and $(b,a)$ has a surjective homomorphism onto the free groupoid
  $\al{F}_{\mathcal{D}}(x,y)$ with two generators in the variety $\mathcal{D}$
  from the classification results of L\'evai and
  P\'alfy~\cite{LevPalfy-binary-min}
  (see Subsection~\ref{subsec-minclone-few-ops}).
\end{enumerate}  
Moreover, all clones in (a), (b), (c), and (m) are minimal.}

\smallskip

The proof makes essential use of the restriction that the binary minimal clones
considered have finite base sets. One of the main ways this finiteness
assumption was exploited in
\cite{Brady-binary-min-taylor, Brady-binary-min-nontaylor} is the following.
Brady used several constructions to produce
better-behaved binary operations in a clone $\CC=\Clo(\al{A}_\circ)$
than the generating operation~$\circ$.
Many of these constructions rely on sequences of
iterated compositions --- like
\mbox{$x, x\circ y, (x\circ y)\circ y, ((x\circ y)\circ y)\circ y,\dots$} ---
which stabilize after
finitely many steps if the base set is finite, and hence 
produce a new operation in $\CC$.
However, if the base set is infinite, then the same conclusions are
not necessarily true.

Concerning the types (B), (S), and (D) in the classification theorem,
which may contain non-minimal clones, Brady also proved in
\cite{Brady-binary-min-taylor, Brady-binary-min-nontaylor} that
if $\mathcal{X}$ is 
the class $\mathcal{BG}_\fin$ of finite Brady groupoids, the class
$\mathcal{SP}_\fin$ of finite spirals, or the class 
$\mathcal{DG}_\fin$ of finite dispersive groupoids, then $\mathcal{X}$
has the following properties: 
\begin{enumerate}
\item[(i)] (Stability under taking subclones)\\
  If $\al{A}_\circ\in\mathcal{X}$ then every subclone of
  $\Clo(\al{A}_\circ)$,
  except the clone of projections, contains a binary operation $\ast$
  such that $\al{A}_\ast=(A;\ast)\in\mathcal{X}$; in particular, every
  minimal subclone of $\Clo(\al{A}_\circ)$ is of the form
  $\Clo(\al{A}_\ast)$ for some $\al{A}_\ast\in\mathcal{X}$.
\item[(ii)] (Smallest witnessing variety)\\
  There exists a finite groupoid $\al{G}_{\mathcal{X}}$ such that
  $\var{V}(\al{G}_{\mathcal{X}})$ is a subvariety of
  $\var{V}(\al{A}_\circ)$ for every $\al{A}_\circ\in\mathcal{X}$.
  The groupoid $\al{G}_{\mathcal{X}}$ can be chosen to be the $2$-generated
  free algebra $\al{F}_{\mathcal{LNB}}(x,y)$, $\al{F}_{\mathcal{SL}}(x,y)$, and
  $\al{F}_{\mathcal{D}}(x,y)$ for $\mathcal{X}=\mathcal{BG}_\fin$,
  $\mathcal{SP}_\fin$, and $\mathcal{DG}_\fin$, respectively.
\item[(iii)] (Membership can be decided efficiently)\\
  There exists a polynomial time algorithm that,
  upon inputting a finite groupoid
  $\al{A}_\circ$, decides whether or not $\al{A}_\circ\in\mathcal{X}$.
\end{enumerate}
In fact, analogous properties~(ii)--(iii)
hold for the classes $\var{V}_\fin$ of finite members of the varieties
$\var{V}=\var{A}_p$ and $p$-$\mathcal{CG}$ (for each prime $p$ separately)
determining the clones in~(a) and~(c), and $\var{V}=\mathcal{RB}$ and
$\var{M}$ determining the clones in~(b)
and~(m). In these cases property~(i) holds trivially,
because the clones are minimal. 

Property~(i) above is the crucial new feature of Brady's broad classification,
which eliminates the `unpredictability' of the types of minimal subclones
of a clone $\CC$ generated by a binary non-projection $\circ$
on a finite set $A$,
provided the clone $\CC$ is `close enough to being minimal' in the sense that
$\al{A}_\circ=(A;\circ)$ belongs to one of the classes
$\mathcal{BG}_\fin$, $\mathcal{SP}_\fin$, or $\mathcal{DG}_\fin$
(which is `easy' to check by condition~(iii)).
In \cite{Brady-binary-min-taylor, Brady-binary-min-nontaylor}, Brady also
investigated the structure of the groupoids in the classes
$\mathcal{BG}_\fin$, $\mathcal{SP}_\fin$, and $\mathcal{DG}_\fin$.

Brady's property~(ii) suffices to imply that the seven families of
clones in the classification are pairwise disjoint, and also helps
to compare the earlier classification theorems for binary minimal clones
with Brady's coarse classification.
For example, if we look at the
classification of minimal clones
$\CC=\Clo(\al{A}_\circ)$ 
with $r:=|\CC^{(2)}|=3$ or~$4$ by
L\'evai and P\'alfy~\cite{LevPalfy-binary-min}
(see Subsection~\ref{subsec-minclone-few-ops}),
we see that
for the two varieties in $\mathcal{F}_3$, the clone 
of a finite groupoid $\al{A}_\circ$ with non-projection operation $\circ$ is
of type~(a) for $\al{A}_\circ$ in $\var{A}_3$ and of type (S) for
$\al{A}_\circ$ in $2$-$\mathcal{SL}$.
For the five varieties in $\mathcal{F}_4$, the clone 
of a finite groupoid $\al{A}_\circ$ with non-projection operation $\circ$ is
of type (c), (b), (B), (D), and (m) for
$\al{A}_\circ$ in $2$-$\mathcal{CG}$, $\mathcal{RB}$,
$2$-$\mathcal{LNB}$, $\var{D}$,
and $\var{M}$, respectively.
The conservative minimal clones classified by
Cs\'ak\'any~\cite{Csakany-conservative-minclone}
(see Subsection~\ref{subsec-minclones-conservative})
are all of the form $\Clo(\al{A}_\circ)$ for a groupoid in the variety
$2$-$\mathcal{SL}$ or $2$-$\mathcal{LNB}$, so --- if the base set is finite ---
they are of type~(S) or~(D).
Among the commutative binary minimal clones classified in
\cite{KSz-commutative-min}
(see Subsection~\ref{subsec-minclone-alg-restr}),
those on a finite base set belong to Brady's type (a), (b), or (c)
in case the clones have abelian representations~\cite{Kearnes-minclone},
and otherwise they are of the form $\Clo(\al{A}_\circ)$ for a groupoid
in the variety $\var{D}$,
$2$-$\mathcal{SL}$ or $2$-$\mathcal{LNB}$, so belong to type (D), (S), or (B).

\section{Minimal clones in other contexts}

Rosenberg, jointly with Machida, made significant contributions
to the study of `essentially minimal'
clones, which are close relatives of minimal clones.
Along with some of their results we will also
discuss another variation on the notion
of minimality for clones,
which is motivated by the Dichotomy Theorem for Constraint
Satisfaction problems (see Subsection~\ref{subsec-csp}), and
a classification of minimal partial clones.

\subsection{Essentially minimal clones}\label{subsec-ess-min}
The minimal clones on a given set $A$ can be thought of as the minimal
members of the family of all `interesting' clones on $A$, provided the clone
$\II_A$ of projections is the only clone on $A$ that is considered
`uninteresting'. The concept of an essentially minimal clone is obtained
the same way if the family of `uninteresting clones' is enlarged to consist
of all essentially unary clones on $A$. Here, a clone is called
\emph{essentially unary}
if all operations in it are essentially unary, and
an operation is essentially unary (see Subsection~\ref{subsec-background})
if it depends on at most one of its
variables. An operation is called \emph{essential} if it depends on at least
two of its variables, i.e., if it is not essentially unary.

So, more formally, the \emph{essentially minimal clones} on a set $A$ can be
defined as the clones $\CC$ on $A$ such that $\CC$ contains
an essential operation, but all proper subclones of $\CC$ are
essentially unary.
It is clear from this definition that a clone $\CC$ is essentially
minimal if and only if $\CC$ is generated by every essential operation in it.
Another immediate consequence of the definition is that
there is a large overlap between minimal and essentially
minimal clones: all minimal clones, except the unary minimal clones,
are also essentially minimal, and an essentially minimal clone is minimal
if and only if it is an idempotent clone. 
Therefore, to separate the study of essentially minimal clones
from the study of minimal clones, the focus in the
study of essentially minimal clones has been on the
non-idempotent essentially minimal clones\footnote{%
In the literature, an ``essentially minimal clone'' is often defined to be
what we call here a non-idempotent essentially minimal clone, and in some papers
(e.g.\ \cite{Machida-Ros-2013}) the name
``essentially minimal clone in the broad sense'' is used for the
clones that we call essentially minimal.},
that is, on the essentially minimal clones generated by
non-idempotent essential operations.

The concept of essentially minimal clones was introduced and first studied
by Machida in \cite{Machida-1981}, and was party motivated by the hope
that the non-idempotent essentially minimal clones are easier to
classify than the minimal clones, due to the presence of a
nontrivial unary operation in the clone.

In \cite{Machida-Ros-1984}, Machida and Rosenberg split the
non-idempotent essentially minimal clones on finite sets $A$
into two types, (A) and (B), and fully classified those
of type (A). To state the classification we will use the following notation.
For an $n$-ary operation $f$ on $A$, $f^*$ denotes the unary operation
$f^*(x):=f(x,\dots,x)$, $\Gamma(f)$ denotes the largest subset of $A$ to which
$f^*$ restricts as a permutation, and $f|_{\Gamma(f)}$ denotes the $n$-ary
function $(\Gamma(f))^n\to A$ obtained from $f$ by restricting the domain to
$(\Gamma(f))^n$. Note that since $A$ is finite, the set $\Gamma(f)$
is nonempty.
A non-idempotent essentially minimal clone $\CC$ is
defined to be of type~(A)
if $\CC$ has a generator $f$ such that the function
$f|_{\Gamma(f)}$ depends on at least two of its variables; otherwise, the type
of $\CC$ is~(B).
The following theorem of Machida and Rosenberg~\cite{Machida-Ros-1984}
classifies the non-idempotent essentially minimal clones of type~(A)
on finite sets. In the theorem statement a clone $\CC$ generated by a
single operation $f$ is called a
\emph{lazy clone} if $\CC$ has no other members than
the projections and the operations obtained from $f$ by
variable manipulations (i.e., by permuting variables, identifying variables,
and adding fictitious variables).

\smallskip
\noindent
{\it 
If $f$ is a non-idempotent $n$-ary operation on a finite set $A$ such that
the function $f|_{\Gamma(f)}$ depends on at least two variables, then
the clone $\CC$ generated by $f$ is essentially minimal
if and only if $f$ satisfies the following conditions:
\begin{enumerate}
\item[{\rm(i)}]
  the identity
  $f(x_1,\dots,x_n)=f(f^*(x_1),\dots,f^*(x_n))$ holds for $f$, and
\item[{\rm(ii)}]
  either
  \begin{enumerate}
  \item[{\rm(a)}]
    the range of $f|_{\Gamma(f)}$ is contained in $\Gamma(f)$ (i.e.,
    $f|_{\Gamma(f)}$ is an operation on $\Gamma(f)$) and
    $f|_{\Gamma(f)}$ generates a minimal clone on $\Gamma(f)$;
  \end{enumerate}  
  or
  \begin{enumerate}  
  \item[{\rm(b)}]
    the range of $f|_{\Gamma(f)}$ is not contained in $\Gamma(f)$, the $n$-ary
    operation\break $f^*(f(x_1,\dots,x_n))$ on $A$ is essentially unary,
    $\CC$ is a \emph{lazy clone}, and $\CC$ is minimal
    among the nontrivial lazy clones on $A$.
  \end{enumerate}
\end{enumerate}  
}
Notice that the identity in~(i) implies that $f^*$ satisfies the
identity $f^*(f^*(x))=f^*(x)$, and since $f$ is not idempotent (that is,
$f^*$ is not the identity map on $A$),
we also have that $\Gamma(f)\subsetneq A$,
$\Gamma(f)$ is the range of $f^*$, and $f^*|_{\Gamma(f)}$ is the identity
map on $\Gamma(f)$.

The non-idempotent essentially minimal clones of type~(B) on finite sets are
not yet classified, but in \cite{Machida-Ros-1993}, Machida and Rosenberg
classified all such clones with binary generators.
In \cite{Machida-Ros-2013} and \cite{Machida-Ros-2014} they also
proved that every non-idempotent essentially minimal clone on a finite set
$A$ is generated by an operation of arity at most $|A|$,
and found all non-idempotent essentially
minimal clones on $3$-element sets, up to similarity.

\subsection{Taylor-minimal clones}\label{subsec-taylor-min}
This notion of minimality for clones 
is motivated by the Dichotomy Theorem for constraint
satisfaction problems $\CSP(\str{A})$, which we discussed in
Subsection~\ref{subsec-csp}. Therefore, we will only consider clones on
finite sets, and will use the terminology and notation introduced in
Subsection~\ref{subsec-csp} without further reference.
In the Dichotomy Theorem --- which proves the Algebraic Dichotomy Conjecture
stated in \eqref{eq-AlgDichConj} for the problems
$\CSP(\str{A})$ where $\str{A}=(A;\rho_1,\dots,\rho_m)$ is a 
finite relational structure with at least two elements and is a core ---
the dividing line between the problems $\CSP(\str{A})$
that are in P and the problems $\CSP(\str{A})$ that are NP-complete is
whether or not the clone $\Pol(\str{A})$ contains a Taylor operation.
The difficult case of the proof is showing that $\CSP(\str{A})$ is in P
if $\Pol(\str{A})$ contains a Taylor operation. Since the complexity
of $\CSP(\str{A})$ can only increase if 
$\Pol(\str{A})$ gets smaller, the clones on $A$ that are minimal for the
property that they contain a Taylor operation play a crucial role.

We will call a clone $\CC$ on a set $A$ a \emph{Taylor clone}
if $\CC$ contains a Taylor operation, and a 
\emph{Taylor-minimal clone} if
$\CC$ is a Taylor clone, but no proper subclone of $\CC$ is a Taylor clone.
In other words, $\CC$ is Taylor-minimal if and only if $\CC$ is Taylor
and is generated by any Taylor operation in it.

Now let $A$ be a finite set.
It is not obvious from the definition of Taylor clones on $A$
--- or from their equivalent
characterization 
(see Subsection~\ref{subsec-csp})
as the clones on $A$ that contain a special weak near unanimity operation
(briefly: special WNU operation) of some arity ---
that every Taylor clone on $A$ has a Taylor-minimal subclone.
However, this fact is an easy consequence of another characterization of
Taylor clones on finite sets, due to Kearnes, Markovi\'{c} and
McKenzie~\cite{Kea-Mar-McK}:

\medskip
\noindent
{\it
  A clone $\CC$ on a finite set $A$ is Taylor if and only if $\CC$ contains
  a $4$-ary idempotent operation $t$ satisfying the
  following ``rare area'' identity in the variables $a,e,r$:
  \begin{equation*}
    t(r,a,r,e)=t(a,r,e,a).
  \end{equation*}%
}%
\noindent
(Such a $t$ is a Taylor operation, so the sufficiency of the condition
is clear.)

A systematic study of Taylor-minimal clones was started very 
recently, and focused so far mainly on how Taylor-minimal clones
can be used to simplify
and unify the two algebraic theories developed by 
Bulatov~\cite{Bulatov} and by
Zhuk~\cite{Zhuk-dich-conj_conf, Zhuk-dich-conj_full}
for their proofs of the Dichotomy Theorem, and a
third theory concerning absorption, which was developed
earlier by Barto and Kozik~\cite{Barto-Kozik1, Barto-Kozik2}
and played a crucial role in understanding which
problems $\CSP(\str{A})$ can be solved by a specific type of
polynomial time algorithm
based on local consistency checking.
All three algebraic theories mentioned here focus on
finite idempotent algebras $\al{A}_\CC:=(A;\CC)$ where
$\CC=\Clo(\al{A})$ is 
a Taylor clone on $A$, and develop results about the clones
$\CC$ themselves 
and/or about the subalgebras of finite powers of $\al{A}$.
Given a coinitial family in the set of all Taylor clones on $A$,\footnote{%
In Zhuk's algorithm that we discussed in Subsection~\ref{subsec-csp},
the family consists of the clones of all algebras
$\al{A}=(A;u)$ where $u$ is a special
WNU operation on $A$.}
the hard part of the Dichotomy Theorem can be proved by finding a
polynomial time algorithm for solving $\CSP(\str{A})$ for
each structure $\str{A}=\str{A}_{\CC,k}$
whose relations are all relations of arity $\le k$
preserving the operations in $\CC$ where $\CC$ is a clone
in the coinitial family and $k\ge2$.
Beyond these common basics,
Bulatov's and Zhuk's algebraic theories don't seem to have much in common;
for example, they use very different coinitial families of Taylor clones,
and while Zhuk's theory does use some concepts from absorption
theory, Bulatov's does not.
Nevertheless, Barto, Brady, Bulatov, Kozik, and Zhuk
discovered in \cite{Barto-Brady-Bulatov-Kozik-Zhuk-confproc,
  Barto-Brady-Bulatov-Kozik-Zhuk}
that by choosing the coinitial family of clones to be the family of
Taylor-minimal clones, a number of connections between the three
theories become apparent. In fact, they say
in~\cite{Barto-Brady-Bulatov-Kozik-Zhuk-confproc}:
``The authors find the extent, to which the notions of the three
theories simplify and unify in minimal Taylor algebras, truly
striking.''

We have seen two characterizations of Taylor clones on finite sets:
one by the existence of a (special) WNU operation in the clone, and the other
by the existence of an idempotent ``rare area'' operation.
There are a number of other, similar characterizations of Taylor clones 
on finite sets, by the existence of a special kind of Taylor operation in
the clone
(e.g., Siggers' $6$-ary operation, see \cite{Siggers}, or a cyclic operation,
see Barto--Kozik~\cite{Barto-Kozik1}).
In \cite{Olsak}, Ol\v{s}\'ak proved a similar characterization for Taylor
clones on arbitrary (possibly infinite) sets:  

\medskip
\noindent
{\it
  A clone $\CC$ on an arbitrary set $A$ is Taylor if and only if $\CC$ contains
  a $6$-ary idempotent operation $f$ satisfying the identities
  \begin{equation*}
    f(x,y,y,y,x,x)=f(y,x,y,x,y,x)=f(y,y,x,x,x,y).
  \end{equation*}%
}%
It follows that Taylor-minimal clones on finite sets contain each of
these types of
Taylor operations. The idempotent ``rare area'' operations 
are optimal among these in terms of their arity, because it was proved
in~\cite{Kea-Mar-McK} that Taylor clones on finite sets
cannot be characterized by the existence of a ternary Taylor operation
satisfying a fixed set of Taylor identities.
In spite of this fact,
Barto, Brady, Bulatov, Kozik, and Zhuk were able to prove
in~\cite{Barto-Brady-Bulatov-Kozik-Zhuk-confproc,
  Barto-Brady-Bulatov-Kozik-Zhuk}
that every Taylor-minimal clone on a finite set $A$ is generated
by a ternary operation.

By inspecting Post's lattice~\cite{Post}, it is easy to see that
on a $2$-element set the Taylor-minimal clones are exactly the non-unary
minimal clones. 
In \cite[Sec.~4.4]{Brady-notes-on-CSP}, Brady classified all conservative
Taylor-minimal clones on finite sets and gave a formula for their number.
He then used this result to classify 
all Taylor-minimal clones on a $3$-element set.
He found that, up to similarity, there are
$24$ Taylor-minimal clones on a $3$-element set, and $19$ of them are
conservative clones.

On finite sets of size $>2$, there exist
Taylor-minimal clones that are not minimal clones.
In fact, there is no upper bound on the height of
the subclone lattices of Taylor-minimal clones on finite sets. For example, 
the \emph{switching operation} $s$ on a set $A$, which is defined to be
the minority operation on $A$ such that $s(a,b,c)=a$
whenever $a,b,c\in A$ are distinct,
generates a Taylor-minimal clone $\CC$ on $A$ whose subclone lattice is a chain
of height $|A|$ if $A$ is finite (a descending $\omega$-chain
if $A$ is infinite), and the unique maximal subclone of $\CC$
is the clone of projections if $|A|=2$, and is generated by a
ternary semiprojection if $|A|>2$
(see~\cite{Marchenkov-homog, szendreiBOOK}).

Nevertheless,
there is some overlap between minimal clones and Taylor-minimal clones.
Clearly, a minimal clone is Taylor-minimal if and only if it is a Taylor clone.
The results we discussed earlier about minimal clones yield that
a minimal clone $\CC$ on a finite set is Taylor-minimal if and only if
(i)~$\CC$ is the clone of an affine module over $\al{Z}_p$ for some prime $p$,
or
(ii)~$\CC$ is a minimal clone in class~(S) of Brady's coarse classification of
binary minimal clones, or
(iii)~$\CC$ is a minimal clone generated by a majority
operation.
Indeed, it is easy to see that the unary minimal clones and
the minimal clones generated by semiprojections
(i.e., the minimal clones of types ${}^\checkmark$(I) and (III))
are not Taylor clones, while the majority minimal clones and the minority
minimal clones (i.e., the minimal clones of types~(IV) and ${}^\checkmark$(V))
are Taylor clones.
Therefore, the only question that remains is this:
Which binary minimal clones contain Taylor operations?
Brady's coarse classification gives a (coarse) answer to
this question. In fact, Brady's coarse classification of binary minimal clones
started in~\cite{Brady-binary-min-taylor}
with a broader goal, namely to understand
the structure of the finite algebras $\al{A}=(A;\CC)$ whose clones
$\CC=\Clo(\al{A})$ are minimal for the property that the corresponding
Constraint Satisfaction Problems $\CSP(\str{A}_{\CC,k})$ (as described earlier
in this subsection) can be solved by local consistency checking.
This property of $\CC$
is closely related to --- though somewhat stronger than ---
the requirement for $\CC$ to be a Taylor-minimal clone.

\subsection{Minimal clones of partial
  operations}\label{subsec-min-partial-clones}
Recall from Subsection~\ref{subsec-max-partial-clones} that
a set $\UU$ of partial operations on a fixed base set $A$ is called
a \emph{clone of partial operations} on $A$, or briefly a
\emph{partial clone} on $A$ if $\UU$
is closed under composition and contains
the (total) projection operations.
Analogously to the definition of a minimal clone (of total operations),
a partial clone $\UU$ on $A$ is called \emph{minimal} if
$\UU$ has exactly one proper partial subclone: the (partial)
clone $\II_A$ of projections.
It follows from these definitions
that every clone of (total) operations is also a partial
clone. In particular, every minimal clone of total operations is a
minimal partial clone.

In \cite{Borner-Haddad-Poschel}, B\"orner, Haddad, and P\"oschel proved that
Statement~\eqref{eq-fin-many-min-clones} extends to partial clones, that is,
on a finite base set $A$, every partial clone other than the
clone $\II_A$ of projections contains a minimal
partial clone, and the number of minimal partial clones is finite.
They also reduced the classification of all minimal partial clones on
finite sets to the classification of all minimal clones of total operations
by proving the following theorem:

\smallskip
\noindent
{\it
  A partial clone $\UU$ on a finite set $A$ is minimal if and only if either
  $\UU$ is a minimal clone of total operations, or $\UU$ is generated
  by a partial projection operation $\pr_i^{(n)}|_D$ where
  $n\ge1$, $1\le i\le n$, and $D\subsetneq A^n$ is an $n$-ary
  totally reflexive, totally symmetric relation on $A$.
}

\bibliographystyle{plain}

\end{document}